\documentclass[11pt]{amsart}
\usepackage{amsmath,amsfonts,amssymb,amsxtra,latexsym,amscd,enumerate,amsthm}

\usepackage{latexsym}
\usepackage{epsfig}
\usepackage{verbatim}
\def\z{\zeta}
\def\t{\theta}
\def\r{\rho}
\def\g{\gamma}
\def\G{\Gamma}
\def\a{\alpha}

\def\d{\delta}
\def\b{\beta}

\def\v{\varphi}
\def\va{\vartheta}
\def\ep{\epsilon}
\def\u{\upsilon}

\def\l{\Lambda}
\def\o{\omega}

\def\n{\nu}
\def\R{\mathbb{R}}
\def\Q{\mathcal{Q}}
\def\C{\mathbb{C}}
\def\Z{\mathbb{Z}}

\def\t{\mathcal{T}}

\def\V{\mathcal{V}}

\def\M{\mathcal{M}}
\def\H{\mathcal{H}}
\def\A{\mathcal{A}}
\def\B{\mathcal{B}}
\def\TT{\mathbb{T}}
\def\N{\mathbb{N}}
\def\n{\mathcal{N}}
\def\f{\mathcal{F}}
\def\Z{\mathbb{Z}}
\def\beq{\begin{equation}}
\def\eeq{\end{equation}}

\def\beq{\begin{equation}}
\def\eeq{\end{equation}}

\newtheorem{t1}{Theorem}
\newtheorem{l1}{Lemma}
\newtheorem{p1}{Proposition}

\begin{document}
\title[]{On the Boundedness of The Bilinear Hilbert Transform along ``non-flat" smooth curves}

\author{Victor Lie}

\date{\today}
\address{Department of Mathematics, Princeton, NJ 08544-1000 USA}

\email{vlie@math.princeton.edu}

\address{Institute of Mathematics of the
Romanian Academy, Bucharest, RO 70700, P.O. Box 1-764, Romania.}

\address{Current address: Department of Mathematics, Purdue University, IN 47907 USA}

\email{vlie@purdue.edu}

\thanks{The author was supported by NSF grant DMS-1200932.}

\maketitle

\begin{abstract}
We are proving $L^2(\R)\times L^2(\R)\,\rightarrow\,L^1(\R)$ bounds for the bilinear Hilbert transform $H_{\G}$ along curves $\Gamma=(t,-\gamma(t))$
with $\g$ being a smooth ``non-flat" curve near zero and infinity.
\end{abstract}
$\newline$

\section{\bf Introduction}

The result that we present here treats the problem of providing $L^2(\R)\times L^2(\R)\,\rightarrow\,L^1(\R)$ bounds for the bilinear Hilbert transform $H_{\G}$ along curves $\Gamma=(t,-\gamma(t))$ defined by
$$H_{\G}: S(\R) \times S(\R)\longmapsto S'(\R)$$
$$H_{\G}(f,g)(x):= \textrm{p.v.}\int_{\R}f(x-t)g(x+\g(t))\frac{dt}{t}\:,$$
where, heuristically, $\g$ is a locally differentiable function which is ``non-flat" (or equivalently,
not ``resembling" a line) near the origin and near the infinity. For the precise statement of the result one should see the next section.

Li (\cite{Li}) was the first to address this topic in the particular case when the curve is given by a monomial $\g(t)=t^d,\,d\in\N\: (d\geq2)$. In this paper we improve his result both \textit{quantitatively}, by obtaining a scale type decay depending on the level sets of the multiplier's phase, and \textit{qualitatively}, by extending the class of functions for which the main theorem holds. Also, our proof relies on different techniques and does not involve the notion of $\sigma-$uniformity used in the monomial case (for an antithesis between the two approaches one should see the Appendix). Instead, and this constitutes one of the novelties in this paper, we will use a discretization procedure of our operator which simultaneously realizes the following: on the one hand separates the variables on the frequency side, on the other hand it preserves the main characteristics (high oscillation and smoothness in one of the variables) of the phase function of the multiplier. This discretization realizes the fragile equilibrium between the two possible extremes: cutting too rough the multiplier (or not at all) which runs into the difficulty of taking advantage of the cancelation offered by the phase (this is one of the reasons for which Li's argument required the subtle concept of $\sigma$
uniformity) or the other scenario of a very fine discretization which comes, among others, with delicate number theoretical problems involving Van der Corput lemma and Weil type sums (also see the discussion in the Appendix).

Another interesting aspect of our proof relies on the class of curves $\gamma$ for which our main theorem holds (see \eqref{asymptotic0}, \eqref{asymptotic0gaminv} and \eqref{fstterm0}). The conditions imposed
on our curves (especially \eqref{asymptotic0} and \eqref{asymptotic0gaminv}) appear as natural but do not seem to have a direct correspondent in the previous math literature (however might be still useful to compare our conditions with the ones found in \cite{NVWW}).

The fundamental concept that comes in the proof of our result is that of \emph{curvature}. It will be responsible for both the strategy chosen to discretize the multiplier of our operator and for the novel - to this problem - ``scale-type decay" (determined by the level sets of the phase function of the multiplier) that we will be able to achieve in order to sum up all the discretized pieces. We mention that our discretization procedure has two stages: the first one - isolating the main term (contribution) of the multiplier - is essentially the only way in which one can take advantage of the multiplier's phase oscillation and thus, it is present in both Li's proof (in a ``more disguised" form) and our approach. The second one, introduced in our paper, relies on a ``separation of variable" argument that will provide the key for proving our Theorem.

This curvature concept proves essential in a series of problems in analysis of which, possibly the most prominent example is offered by the Restriction Problem (\cite{S}). Several of these problems are described in the nice survey paper of Stein and Wainger (\cite{SW}). In particular,
they mention there one classical problem close in spirit to ours, namely the study of the boundedness properties of the Hilbert transform along curves $\Gamma:\:\R\,\rightarrow\,\R^n$, $(n\geq 1)$ defined by
$$\H_{\G}:\: S(\R^n)\,\rightarrow\,S'(\R^n)$$
$$\H_{\G}(f)(x):=p.v. \int_{\R} f(x-\G(t))\,\frac{dt}{t}\,.$$
This topic, as described in \cite{SW}, was initiated by Jones (\cite{Jo}) and  Fabes and Riviere (\cite{FR}) for studying the behavior of the constant coefficient parabolic differential operators. Later it was extended to more general classes of curves (\cite{SW}, \cite{CNSW}) and to
the setting of homogeneous nilpotent Lie groups (\cite{C1}, \cite{C2}).

In a different direction, (see also \cite{Li}), our problem does present some parallelism with the problem of the ($L^2$-)norm-convergence of the nonconventional bilinear averages
$$\frac{1}{N}\,\sum_{n=1}^N f(T^n)\,g(T^{n^2})$$
as $N$ tends to infinity. Here $T$ is an invertible and measure-preserving transformation of a finite measure space $(X,\f,\mu)$ and $f,g\in L^2(X,\f,\mu)$.

As it stands, this topic was treated by Furtstenberg in \cite{F}, and later generalized to powers represented by polynomials (see e.g. \cite{HK}).

However, because of the unsatisfactory  transference principle from the discrete to the continuous setting it does not seem that one can establish
a clear connection between our problem and the corresponding ergodic-theoretical one.

Finally, the study of the boundedness of the bilinear Hilbert transform along curves can be regarded as a natural extension of the work of M. Lacey
and C. Thiele on the resolution of Calderon conjecture (\cite{LT1},\;\cite{LT2}).

Indeed, their situation represents the case of the infinitely
\emph{flat} curve $\g(t)=t$, the absence of the curvature increasing the difficulty level of the proof. For more details as well as other comments relating our problem the reader should consult the Remarks section.
$\newline$

\noindent\textbf{Notations}: In what follows we will repeatedly make use of the following notations. For $A,\,B>0$, we say that $A\lesssim B$ if there exists an absolute constant $C>0$ such that $A\leq C\,B$. In several situations, we choose to write $A\lesssim_{d} B$ in order to stress the $d-$dependence of our constant $C=C(d)>0$ that realizes the inequality $A\leq C\,B$. We will write  $A\approx B$ for two different scenarios which will be clear from the context: either to say that $A$ and $B$ have the same order of magnitude in the sense that $A\lesssim B$ and $B\lesssim A$ or that the size of $A-B$ is much smaller relative to the size of $A$ (the precise quantification of what ``small" means will be specified when needed).
$\newline$

\noindent\textbf{Acknowledgements}: I would like to thank Xiaochun Li and Christoph Thiele for reading the manuscript and providing useful observations. Also special thanks to the referee whose comments improved the presentation of this paper.

\section{\bf Main Result}
$\newline$

In order to state our main result we will first introduce the set $\n\f_{0}$ - the class of all curves $\g$ which are smooth non-flat functions near the origin and obey the following properties:

\begin{itemize}

\item \textit{smoothness, no critical points, variation (near origin)}
\beq\label{nocrtic0}
\eeq
$\exists\: \d>0$ (possibly depending on $\g$) and $\V_{\d}:=(-\d,\d)\setminus\{0\}$ such that $\g\in C^N(\V_{\d})$ $(N\geq4)$ and $|\g'|>0$ on $\V_{\d}$;
moreover
\beq\label{variation0}
\sup_{\a\in\R_{+}}\#\{j\in\Z_{+}\,|\,|2^{-j}\,\g'(2^{-j})|\in[\a,2\a]\}<\infty\:,
\eeq
where here $\Z_{+}:=\{j\in\Z\,|\,j\geq 0\}$.

\item \textit{asymptotic behavior (near origin)}

There exists $\{a_j\}_{j\in\N}\subset\R_{+}$ with $\lim_{j\rightarrow\infty}a_j=0$ such that:
$\newline$

For any $t\in I:=\{s\,|\,\frac{1}{4}\leq|s|\leq 4\}$ and $j\in \Z_{+}$ we have
\beq\label{asymptotic0}
\frac{\g(2^{-j}\,t)}{2^{-j}\,\g'(2^{-j})}=Q(t)\,+\,Q_j(t)\,,
\eeq
with $Q,\,Q_j\in C^N(I)$ and $\|Q_j\|_{C^N(I)}\leq a_j$.

For $s\in J=Q'(I)$ we require
\beq\label{asymptotic0gaminv}
\frac{(\g')^{-1}(s\,\g'(2^{-j}))}{(\g')^{-1}(\g'(2^{-j}))}=r(s)\,+\,r_j(s)\,,
\eeq
where $r,\,r_j\in C^{N-1}(J)$ with $\|r_j\|_{C^{N-1}(J)}\leq a_j\:.$

(The existence of $(\g')^{-1}$, the inverse of $\g'$, will be a consequence of the next hypothesis.)

\item \textit{non-flatness (near origin)}

The main terms in the asymptotic expansion obey
\beq\label{fstterm0}
 \inf_{t\in I} |Q''(t)|,\:\inf_{t\in J} |r'(t)|>c_{\g}>0\:,
\eeq
and
\beq\label{convdualphase0}
 \inf_{{t_1,\,t_2\in J}\atop{t_1\not=t_2}}\frac{|t_1\,r'(t_1)-t_2\,r'(t_2)|}{|t_1-t_2|}>c_{\g}\:.
\eeq
\end{itemize}

In a similar fashion we define $\n\f_{\infty}$ - the class of smooth, non-flat near infinity functions $\g$ having the following properties:

\begin{itemize}

\item \textit{smoothness, no critical points, variation (near infinity)}
\beq\label{nocrticinf}
\eeq
$\exists\: \d>0$ (possibly depending on $\g$) and $\tilde{\V}_{\d}:=(-\infty,-\d)\cup (\d,\infty)$ such that $\g\in C^N(\tilde{\V}_{\d})$ $(N\geq4)$ and $|\g'|>0$ on $\tilde{\V}_{\d}$;
moreover
\beq\label{variationinft}
\sup_{\a\in\R_{+}}\#\{j\in\Z_{-}\,|\,|2^{-j}\,\g'(2^{-j})|\in[\a,2\a]\}<\infty\:,
\eeq
where here $\Z_{-}:=\{j\in\Z\,|\,j\leq 0\}$.

\item \textit{asymptotic behavior (near infinity)}

 There exists $\{\tilde{a}_j\}_{j\in \Z_{-}}\subset\R_{+}$ with $\lim_{j\rightarrow-\infty}\tilde{a}_j=0$ such that:

For any $t\in I:=\{s\,|\,\frac{1}{4}\leq|s|\leq 4\}$ and $j\in \Z_{-}$ we have
\beq\label{asymptoticinfty}
\frac{\g(2^{-j}\,t)}{2^{-j}\,\g'(2^{-j})}=\tilde{Q}(t)\,+\,\tilde{Q}_j(t)\,,
\eeq
with $\tilde{Q},\,\tilde{Q}_j\in C^N(I)$ and $\|\tilde{Q}_j\|_{C^N(I)}\leq \tilde{a}_j$.

For $s\in \tilde{J}=\tilde{Q}'(I)$ we require
\beq\label{asymptoticinftygaminv}
\frac{(\g')^{-1}(s\,\g'(2^{-j}))}{(\g')^{-1}(\g'(2^{-j}))}=\tilde{r}(s)\,+\,\tilde{r}_j(s)\,,
\eeq
where $\tilde{r},\,\tilde{r}_j\in C^{N-1}(\tilde{J})$ with $\|\tilde{r}_j\|_{C^{N-1}(\tilde{J})}\leq \tilde{a}_j\:.$

\item \textit{non-flatness (near infinity)}

The main terms in the asymptotic expansion obey
\beq\label{fstterminfty}
 \inf_{t\in I} |\tilde{Q}''(t)|,\:\inf_{t\in \tilde{J}} |\tilde{r}'(t)|>c_{\g}>0\:,
\eeq
and
\beq\label{convdualphaseinfty}
 \inf_{{t_1,\,t_2\in \tilde{J}}\atop{t_1\not=t_2}}\frac{|t_1\,\tilde{r}'(t_1)-t_2\,\tilde{r}'(t_2)|}{|t_1-t_2|}>c_{\g}\:.
\eeq
\end{itemize}

Finally set $$\n\f:=C(\R\setminus\{0\})\cap\n\f_{0}\cap\n\f_{\infty}$$ and $\n\f^{C}:=\n\f\,+\,Constant$.

The central result in this paper is given by
$\newline$

\noindent\textbf{Main Theorem.} \textit{Let $\G=(t,-\gamma(t))$ be a curve such that $\g\in\n\f^{C}$. Define the bilinear Hilbert transform $H_{\G}$ along the curve $\G$ as $$H_{\G}: S(\R) \times S(\R)\longmapsto S'(\R)$$
$$H_{\G}(f,g)(x):= \textrm{p.v.}\int_{\R}f(x-t)g(x+\g(t))\frac{dt}{t}\:.$$
Then $H_{\G}$ extends boundedly from $L^2(\R)\times L^2(\R)$ to $L^1(\R)$.}
$\newline$

\noindent\textbf{Observations.}

1)  Any real polynomial of degree strictly grater than one belongs to the class $\n\f_{\infty}$ and any real polynomial with no constant and no linear term belongs to $\n\f_{0}$.

2)  If $\g$ is a real analytic function near $0$ (or ${\infty}$) such that $\g(0)=\g'(0)=0$ (or $\g(\infty)=0$) then it belongs to $\n\f_{0}$ ($\n\f_{\infty}$).

3)  Any (real) Laurent polynomial $P(t)=\sum_{j=-n}^m a_j\,t^j$ with $\{a_j\}_j\subset\R$, $a_{-n},\,a_{m}\not=0$ and $n,m\geq 2$ belongs
to $\n\f$.

4)  Any expression (linear combination of terms) of the form $|t|^{\a}\,|\log |t||^{\b}$ with $\a,\,\b\in\R$ and $\a\not\in\{0,\,1\}$ is in $\n\f$.

5)  Notice that $Q':\:I\mapsto Q'(I)$ and $r:\:Q'(I)\mapsto I$ are inverse one to the other and of class $C^{N-1}$.

6)  From the properties of the class $\n\f$ we deduce that
\beq\label{derivativ}
\g\in\n\f\:\:\Rightarrow\:\:\inf_{t\in[\frac{1}{2},\,2]}|Q'(t)|,\:\inf_{t\in[\frac{1}{2},\,2]}|\tilde{Q}'(t)|>0
\eeq
 Moreover, if  $\g\in\n\f$, then there exist $0<C_1(\g)<C_2(\g)$ and $V(0)$ - a neighborhood of $0$ (and respectively infinity - $V(\infty)$) such that $|\g'(\cdot)|>0$ on $V(0)\setminus\{0\}$
(and correspondingly on $V(\infty)$) and
\beq\label{convtr}
 C_1(\g)<\left|\frac{t\,\g''(t)}{\g'(t)}\right|<C_2(\g)\:\:\:\:\:\:\:\:\:\:\:\:\forall\:\:\:t\in (V(0)\setminus\{0\})\cup V(\infty)\:.
\eeq

7) As a consequence of the above observation we have that $\g\in\n\f$ requires the existence of two pairs of constants (depending on $\g$) $K_2\geq K_1>0$ and $C_2\geq C_1>0$ such that for any $t\in V(0)\setminus\{0\}$ (or $t\in V(\infty)$, respectively)
\beq\label{growth}
\eeq
\begin{itemize}
\item either $K_2^{-1}\,|t|^{C_2}<|\g'(t)|<K_1^{-1}\,|t|^{C_1}$;

\item or $\frac{K_1}{|t|^{C_1}}<|\g'(t)|<\frac{K_2}{|t|^{C_2}}$.
\end{itemize}
Deduce that $\g\in\n\f$ implies:

- $\g'$ slowly varying.

- there exist $\lim_{{t\rightarrow 0}\atop{t\not=0}} |\g'(t)|$ and  $\lim_{t\rightarrow \infty} |\g'(t)|$ and both can take only the
values $0$ or $\infty$.

$\newline$
For further discussions see the Remarks section.
$\newline$

Finally, it is worth noticing, that as a consequence of our Main Theorem (see also the second observation above) we have the following:
$\newline$

\noindent\textbf{Corollary.} \textit{Let $\g\in C(\TT)$  be a real function which is real analytic in the origin and such that
$$\g(t)\sim\sum_{k\in N}c_k\,t^k\:\:\:\textrm{when}\:|t|\rightarrow 0,\:\:\:c_{1}=0\:.$$
\indent Then the bilinear Hilbert transform on the torus $H_{\G}^{\TT}$ along the curve $\G=(t,-\gamma(t))$ given by
$$H_{\G}^{\TT}(f,g)(x):= \textrm{p.v.}\int_{\TT}f(x-t)g(x+\g(t))\frac{dt}{t}$$
becomes a bounded operator from $L^2(\TT)\times L^2(\TT)$ to $L^1(\TT)$.}
$\newline$

\noindent\textbf{Observation.} It is known (see \cite{LT2}) that the classical bilinear Hilbert transform $H_{(t,-\a\,t)}$
(here $\a\in\R\setminus\{-1\}$) is a bounded operator from $L^2(\R)\times L^2(\R)$ to $L^1(\R)$.
  If one combines the  Lacey-Thiele proof with that of ours then one can extend the above corollary to the case $c_{1}\in\R\setminus\{-1\}$ (maximal possible range for $c_1$).
$\newline$

\section{\bf The analysis of the multiplier}

If viewed in a multiplier setting, we have:
$$H_{\G}(f,g)(x):= \textrm{p.v.}\,\int_{\R}\int_{\R}\hat{f}(\xi)\hat{g}(\eta)m(\xi,\eta) e^{i\xi x} e^{i\eta x} d\xi d\eta\:.$$ where
$$m(\xi,\eta)=\textrm{p.v.}\,\int_{\R}e^{-i \xi t}\: e^{i \eta \g(t)}\:\frac{dt}{t}\:.$$

Given its singularity and dilation symmetry, we decompose the kernel $\frac{1}{t}$ as follows:
$$\frac{1}{t}=\sum_{j\in\Z} \r_j(t)\:\:\:\:\:\:\:\:\:\forall\:\:t\in\R^{*}\:,$$
where $\r$ is an odd $C^{\infty}$ function with $\operatorname{supp}\,\r\subseteq
\left\{t\in \R\:|\:\frac{1}{4}<|t|<1\right\}$ and $\r_j(t):=2^{j}\r(2^{j}t)$ (with $j\in \Z$).

Consequently,
$$m(\xi,\eta)=\sum_{j\in\Z} m_j(\xi,\eta)$$ with
\beq\label{mj}
m_j(\xi,\eta)=\int_{\R}e^{-i \xi t}\: e^{i \eta \g(t)}\:\r_j(t) dt= \int_{\R}e^{-i\, \frac{\xi}{2^j}\, t}\: e^{i \,\eta \g(\frac{t}{2^j})}\:\r(t)\,dt\:.
\eeq

As one might expect from the above chain of equalities, we need to understand the behavior of $m_j$ when $|j|\rightarrow \infty$ which further involves the properties of $\g$ near the origin and infinity.

In what follows we will only focus on the behavior of $\g$ near the origin, hence from now on, without loss of generality, we assume $j\in\N$.

Our analysis will rely on the essential role played by the curvature condition \eqref{fstterm0}.

As an exemplification for the use of the curvature condition in deriving asymptotic behavior
let us remind the following classical result (see \cite{S} Proposition 3, page 334 and Remark 1.3.4. page 337):
$\newline$

\noindent\textbf{Proposition (model - curvature/asymptotic behavior).}

\noindent \textit{Let $\omega\in C^2((-\d+p,\d+p))$ with $\d,p\in\R$ and $\d>0$ and suppose that $\o'(p)=0$ but $\o''(p)\not=0$.
Then if $\lambda>0$ and $a\in C^{\infty}((-\d+p,\d+p))$ is supported in a \emph{sufficiently small neighborhood of $p$} we have that}
\beq\label{fouriertr}
 \int_{R} e^{-\pi\,i\,\lambda\,\omega(x)}\,a(x)\,dx= c\,e^{-\pi\,i\,\lambda\,\o(p)}\,\lambda^{-\frac{1}{2}}\,|\o''(p)|^{-\frac{1}{2}}\,a(p)
 + O(\lambda^{-\frac{3}{2}})\:.
\eeq
$\newline$
\indent Another instance, of similar flavor, where the curvature condition plays an important role is the following formulation (see e.g. \cite{S}, page 358) for the asymptotic of the Fourier transform of $e^{i \Phi(t)}$ where here $\Phi(t)$ is a $C^1-$real function which, for $t>1$, is strictly convex and obeys $\textrm{lim}_{t\rightarrow\infty}\Phi'(t)=+\infty$:
\beq\label{fouriertrconv}
\int_{1}^{\infty} e^{-i \xi t}\: e^{i \Phi(t) } a(t) dt \approx c\,e^{-i\Psi(\xi)}\,\Phi''(\Psi'(\xi))^{-1/2}\,a(\Psi'(\xi))\:,
\eeq
in the sense that the quotient of both sides tends as $1$ as $\xi\rightarrow \infty$.

Here $\Psi(\xi):=\sup_{t}(t\xi-\Phi(t))$ is the dual phase of $\Phi$, $c$ is an absolute constant, while $a$ is some ``well behaved" smooth function.

\noindent Remark. From the definition one notices that $\Psi'$ and $\Phi'$ are inverses of each other. The concept of duality of phases can be extended to pair of functions that are complementary in the sense of Young's inequality (see \cite{S}, \cite{Z}).

Inspired by \eqref{fouriertr} and  \eqref{fouriertrconv}, our intention is to show that in the setting
$\Phi_{\eta}(t)=\eta\,\g(t)$ and $\Psi_{\eta}$ given by $\Psi'_{\eta}(\cdot)=(\eta\,\g')^{-1}(\cdot)$, for
 $|\Phi_{\eta}(\cdot)|$ large enough, our symbol $m_j$ has the form:
\beq\label{pars}
\int_{\R}e^{-i\, \xi\, t}\: e^{i \,\eta \g(t)}\:\r_j(t)\,dt=
\:c\,e^{-i\Psi_{\eta}(\xi)}\,(\Psi''_{\eta}(\xi))^{1/2}\,2^{j}\,\r^{*}\left(2^j\,(\g')^{-1}(\frac{\xi}{\eta})\right)\,+\,\textrm{Er}\:.
\eeq
where here the error term obeys
$$\textrm{Er}=o((\Psi''_{\eta}(\xi))^{1/2})\;,$$
and $\r^{*}$ is a smooth function, with $\r^{*}\in C^N(\R)$ and $\textrm{supp}\,\r^{*}\subseteq\textrm{supp}\,\r$. (Throughout the paper, the function $\r^{*}$ is allowed to change from line to line.)

Relations \eqref{fouriertr} and \eqref{fouriertrconv} are an expression of the \textit{stationary phase principle} which loosely asserts that the main contribution for an oscillatory integral comes from the information concentrated near the stationary points of the phase.

Following this, when studying our multiplier
$$ m_j(\xi,\eta)=\int_{\R}e^{-i\, \frac{\xi}{2^j}\, t}\: e^{i \,\eta \g(\frac{t}{2^j})}\:\r(t)\,dt\:,$$
we first need to understand the stationary point(s) of the phase function. Thus we are led to the study of the equations of the form $\frac{\xi}{2^j}-\frac{\eta}{2^j}\,\g'(\frac{t}{2^j})=0$ when $|j|$ large and $t\in\textrm{supp}\,\r$. Remark that property \eqref{asymptotic0}, (with the obvious correspondent near infinity) suggests the further analysis of $m_j$ relative to the size of the terms $\frac{\xi}{2^j}$ and $\frac{\eta}{2^j}\,\g'(2^{-j})$.

Based on this we are invited to split $m_j$ as follows:

Let $\nu_0$, $\nu_1$, $\nu_2$ be (even) positive smooth functions such that $\nu_0\in C_{0}^{\infty}(\R)$ with $\textrm{supp}\, \nu_0 \subset (-9/10, 9/10)$, $\nu_1\in C_{0}^{\infty}(\R)$ with $\textrm{supp}\, \nu_1 \subset\{x\,|\,\frac{1}{2}<|x|<2\}$, $\nu_2\in C^{\infty}(\R)$ with $\textrm{supp}\, \nu_2\, \subset \{x\,|\,|x|>3/2\}$ and
$$\nu_0+\nu_1+\nu_2=1\:.$$
Set now $\nu_{j,k}(x):=\nu_k(2^{-j} x)$ and $\nu_{\g,j,k}(x):=\nu_k(2^{-j}\,\g'(2^{-j})\, x),\:k\in \{0,1,2\}$.
(Here we assume wlog that $\|Q\|_{\infty},\:\|\tilde{Q}\|_{\infty}\leq 1$. Otherwise, for the case $j\in\N$, we ``renormalize"
$\nu_{\g,j,k}(x):=\nu_k\large(2^{-j}\,\g'(2^{-j})\,(1+\|Q\|_{\infty})\, x\large),\:k\in \{0,1,2\}$ while for $j\in\Z\setminus\N$ we set
$\nu_{\g,j,k}(x):=\nu_k\large(2^{-j}\,\g'(2^{-j})\,(1+\|\tilde{Q}\|_{\infty})\, x\large),\:k\in \{0,1,2\}$.)

Then, each component $m_j$ of the multiplier $m$ is expressed as:
$$ m_j=\sum_{k,l=0}^{2}\,m_j^{kl}$$
where
$$m_j^{kl}(\xi,\eta):=m_j(\xi,\eta)\, \nu_{j,k}(\xi)\, \nu_{\g,j,l}(\eta)\:.$$

This last relation can be written in a more explicit form:
\beq\label{symbol}
m_j^{kl}(\xi,\eta)=\left(\int_{\R}e^{-i\, \frac{\xi}{2^j}\, t}\: e^{i \,\eta\,\g(\frac{t}{2^j})}\:\r(t)\,dt\right)\, \nu_{k}(\frac{\xi}{2^j})\, \nu_{l}(2^{-j}\,\g'(2^{-j})\,\eta)\:.
\eeq

Using now Taylor expansions and taking advantage of the mean zero property of the function $\r$ together with \eqref{variation0}, \eqref{asymptotic0}, \eqref{variationinft} and \eqref{asymptoticinfty} one notices that for $k,l\in\{0,1\}$ and $j\in\N$ ($j$ large enough) each $m_j^{kl}$ can be essentially reduced to the study of the symbols having the form
\beq\label{reduct}
u_j(\xi,\eta):=\psi(\xi\,2^{-j})\,\varphi(\eta\,2^{-j}\,\g'(2^{-j}))\:,
\eeq
where $\psi,\:\eta$ are smooth compactly supported functions and at least one of them has a zero at the origin (i.e. $\psi(0)=0$ or $\varphi(0)=0$).

Indeed, let us see briefly how can one justify for $k,l\in\{0,1\}$ the reduction of our multiplier to an expression of the form \eqref{reduct}.

First we notice that from the definition of $\nu_0$ and $\nu_1$ we must have that $|\frac{\xi}{2^j}|,\,|\eta\,2^{-j}\,\g'(2^{-j})|< C$.
Moreover, from \eqref{asymptotic0} one also has
$$\sup_{t\in \textrm{supp}\,\r} |\eta\,\g(\frac{t}{2^j})|< C_{\g}\;,$$
where here $C_{\g}>0$ is a constant depending only on $\g$ that is allowed to change from line to line.

Thus it is safe to use Taylor series for writing our symbol as
\beq\label{symbol}
m_j^{kl}(\xi,\eta)=\left(\int_{\R}\sum_{p\in\N} \frac{(-i\, \frac{\xi}{2^j}\, t)^p}{p!}\: \sum_{r\in\N}
\frac{(i \,\eta\,\g(\frac{t}{2^j}))^r}{r!}\:\r(t)\,dt\right)\, \nu_{k}(\frac{\xi}{2^j})\, \nu_{l}(2^{-j}\,\g'(2^{-j})\,\eta)\:.
\eeq
The key point where we use the mean zero condition of $\r$ is when $p=r=0$ to conclude that the corresponding term vanishes.
As a consequence, we deduce that
\beq\label{symbol1}
\eeq
$$m_j^{kl}(\xi,\eta)=\sum_{p+r\geq 1} \frac{(-1)^p\,i^{p+r}}{p!\,r!}\,C_{p,r,j}\:
 \tilde{\nu}_{k,p}(\frac{\xi}{2^j})\: \tilde{\nu}_{l,r}(2^{-j}\,\g'(2^{-j})\,\eta)\,,$$
where
\beq\label{constprj}
C_{p,r,j}:=\int t^p\,\left(\frac{\g(t\,2^{-j})}{2^{-j}\,\g'(2^{-j})}\right)^r\,\r(t)\,dt\:,
\eeq
and $\tilde{\nu}_{k,p}(\xi):=\xi^p\,\nu_{k}(\xi)$, $k\in\{0,\,1\}$.

Since $|C_{p,r,j}|\leq C_{\g}^{p+r}$ the sum in \eqref{symbol1} is absolutely convergent and thus the reduction claimed above holds.

Now taking
$$u=\sum_j u_j\,,$$
(equality above is understood in the sense of distributions) and defining
$$\V(f,g)(x):= \textrm{p.v.}\int_{\R}\int_{\R}\hat{f}(\xi)\hat{g}(\eta)u(\xi,\eta) e^{i\xi x} e^{i\eta x} d\xi d\eta\:,$$
after remodeling the multiplier $u$, one can put together the various techniques used to prove the Coifman-Meyer theorem (\cite{CM1}, \cite{CM2}, \cite{MPTT})
to conclude

\begin{t1}\label{paraproduct}
For any $\frac{1}{p}+\frac{1}{q}=\frac{1}{r}$ with $1<p,q\leq\infty$ and $r\geq 1$ we have
$$\left\|\V(f,g)\right\|_r \lesssim_{\g,p,q,r} \left\|f\right\|_p\left\|g\right\|_q\:.$$
\end{t1}

(We mention here that the in the special case $\g(t)=t^d,\:d\in\N,\:d\geq 2$ the above result can be recovered from \cite{Liunif}.)

Thus, Theorem \ref{paraproduct} solves our concerns relative to the boundedness properties
of the operators given by the symbols
$$m_{kl}:=\sum_{j\in\Z} m_j^{kl}\,,$$
when $k,l\in\{0,1\}$.

The cases $|k-l|\geq1$ and $k=2$ or $l=2$ follow from the decay obtained by applying
the nonstationary phase method (indeed in this case
the phase function appearing in our symbol is highly oscillatory and has no stationary points). This
procedure will be used a bit later when treating the off-diagonal part of the symbol $m_{22}$ (see Claim 1 below)
and thus we will not insist on details here.

Now we turn our attention towards the last possible situation in our decomposition of $m_j$, namely $m_j^{22}$.
Using again ``duality formula" \eqref{pars}, and preserving the previous notations, we claim that
\beq\label{centsymbol}
m_j^{22}(\xi,\eta)=\:c\,e^{-i\Psi_{\eta}(\xi)}\,(\Psi''_{\eta}(\xi))^{1/2}\,2^j\,\r^{*}\left(2^j\,
(\g')^{-1}(\frac{\xi}{\eta})\right)\,+\, \textrm{Error term}\:.
\eeq
(In formula \eqref{centsymbol} we ignore the localization function $\nu_2$ and thus assume that both $\frac{\xi}{2^j}$ and $2^{-j}\,\g'(2^{-j})$ belong to the support of $\nu_2$.)

Indeed, set $\phi$ a smooth compactly supported function with
 \beq\label{suport}
  \textrm{supp}\,\phi\subset\{x\,|\,\frac{1}{10}<|x|<10\}\:.
 \eeq
 Then, using that $\operatorname{supp}\,\r\subseteq \left\{t\in \R\:|\:\frac{1}{4}<|t|<1\right\}$, we rewrite our symbol as:
$$m_j^{22}(\xi,\eta)=\sum_{m,n\in\N} m_{j,m,n }^{22}(\xi,\eta)\,$$
where
\beq\label{piece}
m_{j,m,n}^{22}(\xi,\eta)=\,\left(\int_{\R}e^{-i\, \frac{\xi}{2^j}\, t}\: e^{i \,\eta \g(\frac{t}{2^j})}\:\r(t)\,dt\right)\,\phi\left(\frac{\xi}{2^{m+j}}\right)\,
\phi\left(\frac{\eta}{2^{n+j}}\,\g'(\frac{1}{2^j})\right)\:.
\eeq
(From now on, for notational simplicity, the function $\nu_2$ will be absorbed in the function $\phi$.)

Now based on the intuition built previously, we expect two different regimes: the stationary phase regime - represented by the ``diagonal" terms ($m\approx n$) and the non-stationary regime - terms living far from diagonal, ie $|m-n|>C(\g)>1$.

The terms far from diagonal are only contributing to the error term in \eqref{centsymbol}, and they are treated as follows:
 \begin{itemize}
\item{ \textbf{Claim 1.}  Since $|m-n|>>1$ integrating once by parts we have
 \beq\label{errterm}
m_{j,m,n}^{22}(\xi,\eta)=\frac{1}{2^{\max\{m,n\}}}\,\tilde{m}_{j,m,n}^{22}(\xi,\eta)\:,
\eeq
where $\tilde{m}_{j,m,n}^{22}$ is a symbol defined by \eqref{piece} with just replacing $\r$ and $\phi$ with different same-structured smooth functions. If we set now $\V_{j,m,n}$ the bilinear operator having the symbol given by $m_{j,m,n}^{22}$, based on \eqref{errterm} for each fixed
$j$ we have
\beq\label{partoperat}
\left\|\V_{j,m,n}(f,g)\right\|_1 \lesssim_{\g}\frac{1}{2^{\max\{m,n\}}}\, \left\|\hat{f}(\cdot)\, \phi_1\left(\frac{\cdot}{2^{m+j}}\right)\right\|_2\left\|\hat{g}(\cdot)\,
\phi_2\left(\frac{\cdot}{2^{n+j}}\,\g'(\frac{1}{2^j})\right)\right\|_2\:.
\eeq}
\item{ \textbf{Claim 2.} Taking now $\V_{m,n}(f,g):=\sum_{j}\V_{j,m,n}(f,g)$ and making essential use of hypothesis \eqref{variation0}, \eqref{variationinft} and relation \eqref{partoperat}, by Cauchy-Schwarz inequality we have
\beq\label{sumj}
\left\|\V_{m,n}(f,g)\right\|_1 \lesssim_{\g}\frac{1}{2^{\max\{m,n\}}}\, \|f\|_2\,\|g\|_2\,,
\eeq
which further implies
\beq\label{partoperat1}
\left\|\sum_{{m,n}\atop{|m-n|>>1}}\V_{m,n}(f,g)\right\|_1\lesssim \|f\|_2\,\|g\|_2\,.
\eeq}
\end{itemize}
Indeed, in what follows we provide a concise justification for each of the above claims:

\noindent\textbf{Proof of Claim 1.} Define the differential operator
\beq\label{difop}
L:=\frac{-i}{-\frac{\xi}{2^j}+\frac{\eta}{2^j}\,\g'(\frac{t}{2^j})}\,\partial_{t}\:.
\eeq
This operator has the key property that
\beq\label{fixpoint}
L(e^{-i\, \frac{\xi}{2^j}\, t}\: e^{i \,\eta \g(\frac{t}{2^j})})=e^{-i\, \frac{\xi}{2^j}\, t}\: e^{i \,\eta \g(\frac{t}{2^j})}\:.
\eeq
Making use of \eqref{fixpoint} and integrating by parts we have
\beq\label{partintegr}
\int_{\R} L(e^{-i\, \frac{\xi}{2^j}\, t}\: e^{i \,\eta \g(\frac{t}{2^j})})\,\r(t)\,dt=
\int_{\R} e^{-i\, \frac{\xi}{2^j}\, t}\: e^{i \,\eta \g(\frac{t}{2^j})}\,L^{\tau}(\r(t))\,dt
\eeq
where we set
$$L^{\tau}(\r(t)):= \partial_{t}\left( \frac{i}{-\frac{\xi}{2^j}+\frac{\eta}{2^j}\,\g'(\frac{t}{2^j})}\,\r(t)\right)\:.$$
Define
$$A_{j,m,n}(\xi,\eta):=$$
$$\left(\int_{\R} e^{-i\, \frac{\xi}{2^j}\, t}\: e^{i \,\eta \g(\frac{t}{2^j})}\,\frac{i\,\r'(t)}{-\frac{\xi}{2^j}+\frac{\eta}{2^j}\,\g'(\frac{t}{2^j})}\,dt\right)\,\phi\left(\frac{\xi}{2^{m+j}}\right)\,
\phi\left(\frac{\eta}{2^{n+j}}\,\g'(\frac{1}{2^j})\right)\:, $$
and respectively
$$B_{j,m,n}(\xi,\eta):=$$
$$\left(\int_{\R} e^{-i\, \frac{\xi}{2^j}\, t}\: e^{i \,\eta \g(\frac{t}{2^j})}\,\frac{-i\,\frac{\eta}{2^{2j}}\,\g''(\frac{t}{2^{j}})}{(-\frac{\xi}{2^j}+\frac{\eta}{2^j}\,\g'(\frac{t}{2^j}))^2}\,\r(t)\,dt\right)\,
\phi\left(\frac{\xi}{2^{m+j}}\right)\,
\phi\left(\frac{\eta}{2^{n+j}}\,\g'(\frac{1}{2^j})\right)\:.$$
Thus we have
$$m_{j,m,n}^{22}(\xi,\eta)= A_{j,m,n}(\xi,\eta)\,+\,B_{j,m,n}(\xi,\eta)\;.$$
\textbf{Case 1}. $m>>_{\g} n$ and $j\in\N$ large.
$\newline$

\noindent Remark. We will only consider large values of $j$
for which $\|Q_j\|_{C^N(I)}\leq a_j\leq \frac{1}{2}\inf_{t\in I} |Q'(t)|\:.$

Applying Taylor series (notice that for the series to be convergent we will require $m-n> 1000\,\|Q'\|_{C(I)}$ - see \eqref{asymptotic0}) we have
$$\frac{1}{-\frac{\xi}{2^j}+\frac{\eta}{2^j}\,\g'(\frac{t}{2^j})}=-\frac{1}{2^m}\,\frac{1}{\frac{\xi}{2^{j+m}}}\,\sum_{l=0}^{\infty} \frac{1}{2^{l(m-n)}}\,\left(\frac{\frac{\eta}{2^{n+j}}\,\g'(\frac{1}{2^j})}{\frac{\xi}{2^{j+m}}}\right)^l\,\left(\frac{\g'(\frac{t}{2^j})}{\g'(\frac{1}{2^j})}\right)^l\;,$$
and hence
$$A_{j,m,n}=\sum_{l\in\N} A_{j,m,n,l}\,,$$
where
$$A_{j,m,n,l}(\xi,\eta):=-\frac{1}{2^m}\,\frac{1}{2^{l(m-n)}}$$
$$\left(\int_{\R} e^{-i\, \frac{\xi}{2^j}\, t}\: e^{i \,\eta \g(\frac{t}{2^j})}\,i\,\r'(t)\,(Q'(t)+Q'_j(t))^l\,dt\right)\,
\phi_l\left(\frac{\xi}{2^{m+j}}\right)\,\tilde{\phi}_l\left(\frac{\eta}{2^{n+j}}\,\g'(\frac{1}{2^j})\right)\,$$
with $\phi_l,\,\tilde{\phi}_l$ smooth, compactly supported away from the origin and $\|\phi_l\|_{C^{\a}},\,$
$\|\tilde{\phi}_l\|_{C^{\a}}\lesssim \a!\,|l|^{\a}\,|C|^l$ for some fixed $C\in\R$.
From this we deduce that the bilinear operator
\beq\label{form}
\l^{A}_{j,m,n,l}(f,g)(x):=\int_{\R}\int_{\R} \hat{f}(\xi)\,\hat{g}(\eta)\,A_{j,m,n,l}(\xi,\eta)\,e^{i\,\xi\,x}\,e^{i\,\eta\,x}\,d\xi\,d\eta
\eeq
obeys
$$\l^{A}_{j,m,n,l}(f,g)(x)=-\frac{1}{2^m}\,\frac{1}{2^{l(m-n)}}\,\times$$
$$\int_{\R}\,
\left(\hat{f}(\cdot)\,\phi_l(\frac{\cdot}{2^{m+j}})\right)^{\check{}}\,(x-\frac{t}{2^j})\:\left(\hat{g}(\cdot)\,
\tilde{\phi}_l(\frac{\cdot\,\g'(2^{-j})}{2^{n+j}})\right)^{\check{}}\,(x+\g(\frac{t}{2^j}))\,i\,\r'(t)\,(Q'(t)+Q'_j(t))^l\,dt\:.$$
Thus applying Cauchy-Schwarz followed by Parseval we have
\beq\label{normform}
\|\l^{A}_{j,m,n,l}(f,g)\|_{1}\lesssim \frac{1}{2^m}\,\left(\frac{2\,\|Q'\|_{\infty}}{2^{m-n}}\right)^l\,\|\hat{f}(\cdot)\,\phi_l(\frac{\cdot}{2^{m+j}})\|_2\,
\|\hat{g}(\cdot)\,\tilde{\phi}_l(\frac{\cdot\,\g'(2^{-j})}{2^{n+j}})\|_2\:.
\eeq
A similar argument applies to the multiplier $B_{j,m,n}$ with the extra twist that from \eqref{variation0} one
has that $\frac{1}{2^{2j}}\,\g''(\frac{t}{2^{j}})= 2^{-j}\,\g'(2^{-j})\,(Q''(t)+Q_j''(t))$. With the obvious changes, same reasonings
apply in the case $j\in\Z\setminus \N$.

From the facts described above, we conclude that \eqref{partoperat} holds.
$\newline$
\noindent\textbf{Case 2}. $n>>_{\g}m$ and $j\in\N$ large.
$\newline$

As before, we apply a Taylor series argument. This time, for the convergence of our series, we require $n-m> 1000\,\|\frac{1}{Q'}\|_{C(I)}$. Notice that \eqref{nocrtic0}, \eqref{asymptotic0} and \eqref{fstterm0} assures that  $\|\frac{1}{Q'}\|_{C(I)}<\infty$.

 These being said, we first have
$$\frac{1}{-\frac{\xi}{2^j}+\frac{\eta}{2^j}\,\g'(\frac{t}{2^j})}=\frac{1}{2^n}\,\frac{1}{\frac{\eta}{2^{j+n}}\,\g'(\frac{1}{2^j})}\,
\frac{\g'(2^{-j})}{\g'(t 2^{-j})}\,\sum_{l=0}^{\infty} \frac{1}{2^{l(n-m)}}\,\left(\frac{\frac{\xi}{2^{m+j}}}{\frac{\eta\,\g'(\frac{1}{2^j})}{2^{j+n}}}\right)^l\,
\left(\frac{\g'(\frac{1}{2^{j}})}{\g'(\frac{t}{2^j})}\right)^l\;.$$
In this setting we further have
$$A_{j,m,n}=\sum_{l\in\N} \tilde{A}_{j,m,n,l}\,,$$
where
$$\tilde{A}_{j,m,n,l}(\xi,\eta):=\frac{1}{2^n}\,\frac{1}{2^{l(n-m)}}\times$$
$$\left(\int_{\R} e^{-i\, \frac{\xi}{2^j}\, t}\: e^{i \,\eta \g(\frac{t}{2^j})}\,i\,\frac{\r'(t)}{(Q'(t)+Q'_j(t))^{l+1}}\,dt\right)\,
\phi_l\left(\frac{\xi}{2^{m+j}}\right)\,\tilde{\phi}_l\left(\frac{\eta}{2^{n+j}}\,\g'(\frac{1}{2^j})\right)\,$$
and, as before, $\phi_l,\,\tilde{\phi}_l$ smooth, compactly supported away from the origin and $\|\phi_l\|_{C^{\a}},\,$
$\|\tilde{\phi}_l\|_{C^{\a}}\lesssim \a!\,|l|^{\a}\,|C|^l$.
We therefore have that
\beq\label{form2}
\l^{\tilde{A}}_{j,m,n,l}(f,g)(x):=\int_{\R}\int_{\R} \hat{f}(\xi)\,\hat{g}(\eta)\,\tilde{A}_{j,m,n,l}(\xi,\eta)\,e^{i\,\xi\,x}\,e^{i\,\eta\,x}\,d\xi\,d\eta
\eeq
obeys
$$\l^{\tilde{A}}_{j,m,n,l}(f,g)(x)=\frac{1}{2^n}\,\frac{1}{2^{l(n-m)}}\,\times$$
$$\int_{\R}\,
\left(\hat{f}(\cdot)\,\phi_l(\frac{\cdot}{2^{m+j}})\right)^{\check{}}\,(x-\frac{t}{2^j})\:\left(\hat{g}(\cdot)\,
\tilde{\phi}_l(\frac{\cdot\,\g'(2^{-j})}{2^{n+j}})\right)^{\check{}}\,(x+\g(\frac{t}{2^j}))\,i\,\frac{\r'(t)}{(Q'(t)+Q'_j(t))^{l+1}}\,dt\:.$$
Now since $Q'(t)+Q'_j(t)$ is uniformly bounded away from the origin (based on the hypothesis imposed on the curve $\g$), we can apply the same reasonings as before, i.e Cauchy-Schwarz followed by Parseval, to deduce
\beq\label{normform2}
\|\l^{\tilde{A}}_{j,m,n,l}(f,g)\|_{1}\lesssim \frac{1}{2^n}\,\left(\|\frac{2}{Q'}\|_{\infty}\,\frac{1}{2^{n-m}}\right)^l\,\|\hat{f}(\cdot)\,\phi_l(\frac{\cdot}{2^{m+j}})\|_2\,
\|\hat{g}(\cdot)\,\tilde{\phi}_l(\frac{\cdot\,\g'(2^{-j})}{2^{n+j}})\|_2\:.
\eeq
Summing over the parameter $l$, we obtain a good control of the bilinear operator having $A_{j,m,n}$ as a multiplier. Similar reasonings apply
for controlling the bilinear multiplier represented by $B_{j,m,n}$. All these are then easily extended to the case $j\in\Z\setminus\N$.
$\newline$
\noindent\textbf{Proof of Claim 2.}
For proving \eqref{partoperat1} we only need to show that
\beq\label{multiplier}
\sum_{j\in\Z}\|\hat{g}(\cdot)\,\phi(\frac{\cdot}{2^{n+j}}\,\g'(2^{-j}))\|_2^2\lesssim \|g\|_2^2\,:.
\eeq
This further reduces to showing that
$$\sum_{j\in\Z} |\phi (\frac{\cdot}{2^{n+j}}\,\g'(2^{-j}))|^2\lesssim 1\:,$$
holds uniformly in $n$. Here, as before, $\phi$ is a $C_{0}^{\infty}$ function compactly supported away from the origin.

But this is equivalent with proving that, for some $c>1$ fixed, we have
\beq\label{linf}
 \|\sum_{j\in\Z} \chi_{[\frac{1}{c\,\g'(2^{-j})\,2^{-j-n}},\frac{c}{\g'(2^{-j})\,2^{-j-n}}]}\|_{\infty}\lesssim 1\;.
\eeq
But this is a direct consequence of \eqref{variation0}.

With this we are done with proving the Claims 1 and 2.
$\newline$

We now turn our attention towards the diagonal term.

Thus, for $m,\,n\in\N$ with
\beq\label{mn}
|m-n|\leq C(\g)=1000\max\{\|Q'\|_{C(I)},\,\|\frac{1}{Q'}\|_{C(I)}\}\,,
\eeq
we analyze the term $m_{j,m,n}^{22}(\xi,\,\eta)$.

$\newline$
\textbf{Analyzing the multiplier's phase.}
$\newline$

In this subsection we focus on the main properties of the phase
$$\varphi_{\xi,\eta}(t):=-\frac{\xi}{2^j}\,t+\eta\,\g(\frac{t}{2^j})\:.$$

From the hypothesis imposed on our curve $\g$, we claim that for $\xi,\,\eta,\,j$ fixed ($j$ large, depending on $\g$), there exists exactly one critical point
\beq\label{critpoint}
t_c=t_c(\xi\,,\eta,\,j)\in [2^{-k(\g)}, 2^{k(\g)}]\,\:\;\textrm{such that}\:\;\varphi'_{\xi,\eta}(t_c)=\frac{\xi}{2^j}-\frac{\eta}{2^j}\,\g'(\frac{t_c}{2^j})=0\;,
\eeq
where here $k(\g)\in\N$ is an integer depending only on $\g$ that will be chosen later (see \eqref{Qprop}).

To see this, we first notice that the function $t\,\rightarrow\,\frac{\g'(t\,2^{-j})}{\g'(2^{-j})}$ is strictly monotone on the interval
$[2^{-k(\g)}, 2^{k(\g)}]$. Indeed, suppose $t\in I_k=[2^{k-2},\,2^{k+2}]\subseteq [2^{-k(\g)}, 2^{k(\g)}]$ and $k\in\Z$. Then, writing $t=2^k\,s$ and applying \eqref{asymptotic0} for $s\in I$ we have
\beq\label{monotone}
\frac{\g'(t\,2^{-j})}{\g'(2^{-j})}=\frac{\g'(s\,2^{-j+k})}{\g'(2^{-j+k})}\cdot\frac{\g'(2^{-j+k})}
{\g'(2^{-j})}=[Q'(s)+Q'_{j-k}(s)]\,\frac{\g'(2^{-j+k})}{\g'(2^{-j})}\:.
\eeq
Since from \eqref{asymptotic0} and \eqref{fstterm0} we have that $Q'(s)+Q'_{j-k}(s)$ is strictly monotone on $I$, the same property holds for $\frac{\g'(t\,2^{-j})}{\g'(2^{-j})}$ on $I_k$. Now, varying $k$ between $-k(\g)$ and $k(\g)$ we obtain the strict monotonicity of $\frac{\g'(t\,2^{-j})}{\g'(2^{-j})}$ on the entire interval $[2^{-k(\g)}, 2^{k(\g)}]$.
 This shows that there is \textit{at most} one critical point $t_c$ obeying \eqref{critpoint}.

Now, to show its existence, it is enough to prove that
\beq\label{range}
[2^{m-n-100},\,2^{m-n+100}]\subseteq \textrm{Range}_{t\in [2^{-k(\g)}, 2^{k(\g)}]}(\frac{\g'(t\,2^{-j})}{\g'(2^{-j})})\:.
\eeq
Applying again \eqref{asymptotic0} we have
$$\frac{\g'(2^{-j+1})}{\g'(2^{-j})}=Q'(2)+Q'_j(2)\:\:\:\textrm{and}\:\;\:\frac{\g'(2^{-j-1})}{\g'(2^{-j})}=Q'(1/2)+Q'_j(1/2)\;.$$
Assuming wlog that $Q'$ is strictly increasing, we have that $0<Q'(1/2)<Q'(1)=1<Q'(2)$ and hence if $j$ large enough choosing $k(\g)$
such that both
\beq\label{Qprop}
Q'(2)^{k(\g)}>2^{1000 C(\g)}\:\:\:\:\:\textrm{and}\:\:\:\:\:Q'(\frac{1}{2})^{k(\g)}< 2^{-1000 C(\g)}
\eeq
(with $C(\g)$ as defined in \eqref{mn}) we deduce that \eqref{range} holds.

 Notice that the existence of such a $t_c$ requires
 \beq\label{tc}
 t_c=2^{j}\,(\g')^{-1}\left(\frac{\xi}{\eta}\right)\:\:\textrm{with}\:\: 2^{m-n-100}\,|\g'(2^{-j})|<|\frac{\xi}{\eta}|<2^{m-n+100}\,|\g'(2^{-j})|\:.
 \eeq

 Also, based on \eqref{asymptotic0} and \eqref{fstterm0}, we notice that
 $$\varphi''_{\xi,\eta}(t)=\frac{\eta}{2^{2j}}\,\g''(\frac{t}{2^j})=\frac{\eta}{2^{j}}\,\g'(\frac{t}{2^j})\,k_j(t)\,,$$
 with $k_j\in C^2(\R)$ and $0<c_{\g}<|k_j(t)|<C_{\g}$ for any $t\in 2\,\textrm{supp}\,\r$.

 Thus, using now \eqref{fstterm0}, we deduce that the Hessian obeys
\beq\label{hess}
|\varphi''_{\xi,\eta}(t)|\gtrsim |\frac{\eta}{2^{j}}\,\g'(\frac{t}{2^j})|\gtrsim 2^m\approx_{\g} 2^n\:\:\:\:\:\:\;\;\:\:\:\:\:\:\forall\:(\xi,\,\eta)\in\textrm{supp}\,m_{j,m,n}^{22}(\cdot,\cdot)\;.
\eeq

This ends our discussion about the analysis of the phase $\varphi_{\xi,\eta}$.

$\newline$
\textbf{Main term of the multiplier.}
$\newline$

We now move further, and present \textbf{the heuristic} of our approach.

\noindent Remark. We only consider here the case $m=n$. Also assume wlog that our integration in $t$ is done only over $\R_{+}$.

As expected, we will make use of (non)stationary phase principle and split the region of integration depending on the behavior (oscillation) of the phase $\varphi_{\xi,\eta}(t)$.

Indeed, for an appropriate choice of the set $V$, we decompose
\beq\label{heuristic}
\eeq
$$m_{j,m,m}^{22}(\xi,\,\eta)\approx e^{i\,\varphi_{\xi,\eta}(t_c)}\,\int_{V} e^{i\,\int_{t_c}^{t}\varphi'_{\xi,\eta}(u)\,du}\r(t)\,dt
\,+\,$$
$$\int_{V^c} e^{i\,\varphi_{\xi,\eta}(t)}\,\frac{\r^{*}(t)}{\eta\, 2^{-j}\,\g'(2^{-j})\,\int_{t_c}^{t}\frac{2^{-2j}
\g''(s 2^{-j})}{2^{-j} \g'(2^{-j})}\,ds}\,dt= A(\xi,\eta)\,+\,B(\xi,\eta)\:.$$

\noindent Remark. The scheme decomposition in \eqref{heuristic} should be only understood at the heuristic level. For the precise decomposition one should read the following paragraph.

\textit{Roughly} speaking, we will have to deal with two components:
\begin{itemize}
\item{$A(\xi,\eta)$ - the main term, catches the information of the phase \textbf{close} to the critical point $t_c$. Indeed, we will construct the neighborhood $V$ such that our phase has essentially no oscillation
     $$|\int_{t_c}^{t}\varphi'_{\xi,\eta}(u)\,du|\lesssim 1\:\:\:\:\:\:\:\forall\:\:t\in V\,.$$
     Based on \eqref{asymptotic0}, for $m$ large enough, this is virtually equivalent with imposing
     $$V=V(\xi,\eta):=\{t\,|\,|t-t_c|\lesssim 2^{-\frac{m}{2}}\}\;.$$
 With this choice we will then show that
 \beq\label{A}
A(\xi,\eta)=e^{i\,\varphi_{\xi,\eta}(t_c)}\,2^{-\frac{m}{2}}\,\r^{*} (t_c)\:,
\eeq
with $\r^{*}$ a smooth function with same properties as $\r$.}
 \item{$B(\xi,\eta)$ - can be regarded as an error term, as it encapsulates the behavior of our multiplier \textbf{far} from the stationary point. Thus integrating by parts we expect some extra decay to come into play. Indeed, using \eqref{asymptotic0} and \eqref{fstterm0} one has:
$$B(\xi,\eta)=\frac{1}{\eta 2^{-j} \g'(2^{-j})}\,\int_{V^{c}}  e^{i\,\varphi_{\xi,\eta}(t)}\,
\frac{\r^{*}(t)}{\int_{t_c}^{t} [Q''(s)+Q''_{j}(s) ]\,ds}\,dt\,,$$
which implies
\beq\label{B}
B(\xi,\eta)\approx\frac{1}{2^{\frac{m}{2}}}\int_{|t-t_c|\gtrsim 2^{-\frac{m}{2}}} e^{i\,\varphi_{\xi,\eta}(t)}\r^{*} (t)\,dt\:.
\eeq
Due to the good decay this last term appears as an error term and can be treated as previously the non-diagonal terms.}
\end{itemize}

We now initiate the final stage of our multiplier decomposition by making rigorous the heuristic presented in \eqref{heuristic}.

Choose $\va,\,\tilde{\va}\in C_{0}^{\infty}(\R)$ such that $\textrm{supp}\,\va\subseteq[-10,10]$,  $\textrm{supp}\,\tilde{\va}\subseteq \{t\,|\,\frac{1}{100}<|t|<100\}$ and
$$1=\va(t)\,+\,\sum_{k\in\N}\tilde{\va}(2^{-k}\,t)\:\:\:\:\:\:\:\:\:\:\:\:\forall\:t\in\R\:.$$
Write now
\beq\label{mdecomp}
m_{j,m,m}^{22}(\xi,\eta)= \A_{j,m}(\xi,\eta)\,+\,\sum_{k=0}^{\infty}\B_{j,m}^k(\xi,\eta)\:,
\eeq
where
\beq\label{Aterm}
\A_{j,m}(\xi,\eta)= \,\left(\int_{\R} \va(2^{\frac{m}{2}}(t-t_c))\,e^{i\,\varphi_{\xi,\eta}(t)}\:\r(t)\,dt\right)\,\phi\left(\frac{\xi}{2^{m+j}}\right)\,
\phi\left(\frac{\eta}{2^{m+j}}\,\g'(\frac{1}{2^j})\right)\:,
\eeq
and
\beq\label{Bkterm}
\B_{j,m}^k(\xi,\eta)= \,\left(\int_{\R} \tilde{\va}(2^{\frac{m}{2}-k}(t-t_c))\,e^{i\,\varphi_{\xi,\eta}(t)}\:\r(t)\,dt\right)\,\phi\left(\frac{\xi}{2^{m+j}}\right)\,
\phi\left(\frac{\eta}{2^{m+j}}\,\g'(\frac{1}{2^j})\right)\:.
\eeq
Notice that since $t_c\in [2^{-k(\g)}, 2^{k(\g)}]$ only finitely many terms $\B_{j,m}^k$ are nonzero (those for which $k\lesssim C(\g)+\frac{m}{2}$; as before here $C(\g)$ designates a (possibly large) positive constant depending only on $\g$).

We pass to the analysis of the term $\A_{j,m}$.

First, we notice that
\beq\label{phasedecomp}
e^{i\,\varphi_{\xi,\eta}(t)}= e^{i\,\varphi_{\xi,\eta}(t_c)}\,e^{i\,\int_{t_c}^{t}\,\int_{t_c}^{s} \varphi''_{\xi,\eta}(r)\,dr\,ds}\:.
\eeq
Making now in \eqref{Aterm} the change of variable $t\,\rightarrow\,t_c+ 2^{-\frac{m}{2}}t$, and based on \eqref{phasedecomp},
we deduce that
\beq\label{Arewrit}
\A_{j,m}(\xi,\eta)= \,2^{-\frac{m}{2}}\,e^{i\,\varphi_{\xi,\eta}(t_c)}\,w_{j,m,\eta}(t_c)\,\phi\left(\frac{\xi}{2^{m+j}}\right)\,
\phi\left(\frac{\eta}{2^{m+j}}\,\g'(\frac{1}{2^j})\right)\:,
\eeq
where for $u\in[2^{-k(\g)},\,2^{k(\g)}]$ we set
\beq\label{wdef}
w_{j,m,\eta}(u):=\int_{\R} \va(t)\,
e^{i\,\int_{0}^{2^{-\frac{m}{2}}\,t}\,\int_{0}^{s} \varphi''_{\xi,\eta}(r+u)\,dr\,ds}\:\r(2^{-\frac{m}{2}}t+u)\,dt\;.
\eeq
The function $w_{j,m,\eta}$ is smooth, and behaves well with respect of all the parameters $j,m,\eta$. More precisely,
applying again a Taylor series argument we have
\beq\label{wcontrol}
w_{j,m,\eta}(u)=\sum_{l=0}^{\infty} \frac{1}{l!}\,\left(\frac{i\,\eta\,\g'(2^{-j})}{2^{m+j}}\right)^l\,\t_{j,m,l}(u)\,,
\eeq
with
\beq\label{tau}
\t_{j,m,l}(u):=\int_{\R} \va(t)\, \left(\int_{0}^t\int_{0}^s \frac{\g''(\frac{r 2^{-m/2}+u}{2^j})}{2^j\,\g'(2^{-j})} dr\,ds\right)^l\,\r(2^{-\frac{m}{2}}t+u)\,dt\:. \eeq
The key facts in \eqref{wcontrol} are:

- $\frac{\eta\,\g'(2^{-j})}{2^{m+j}}\in \textrm{supp}\,\phi$;

- $\t_{j,m,l}$ does not depend on $\eta$;

- $\|\t_{j,m,l}\|_{C^{N-3}([2^{-k(\g)},\,2^{k(\g)}])}\leq (C_{\g})^l$.

Next remark that,  one can extend  \eqref{asymptotic0gaminv} for values of  $s$ not necessarily inside of the interval $J$.
Indeed, one can use similar techniques with those used in the (sub)section ``Analyzing the multiplier's phase" to see that
\beq\label{tcinv}
t_c=2^{j}\,(\g')^{-1}(\frac{\xi}{\eta})= \bar{r}(\frac{\xi}{\eta\,\g'(2^{-j})})\,+\,\bar{r_j}(\frac{\xi}{\eta\,\g'(2^{-j})})\:,
\eeq
where here $\bar{r}$ and $\bar{r_j}$ are extensions of $r,\,r_j$ derived from the initial definition of $r$ and $r_j$ and consistent with the property that $r(s)=2\,r(Q'(\frac{1}{2})\,s)$ for any $s\in Q'(I_0)$ with $I_0=\{s\,|\,\frac{1}{2}\leq |s|\leq 2\}$.

Now, setting $\phi^l(x):=x^l\,\phi(x)$ and combining the above elements, we have:
\beq\label{w}
\A_{j,m}(\xi,\eta)=\sum_{l\in\N} \,\frac{i^l}{l!}\,\A_{j,m,l}(\xi,\eta)\,,
\eeq
with
\beq\label{aa}
\A_{j,m,l}(\xi,\eta):=\,2^{-\frac{m}{2}}\,e^{i\,\varphi_{\xi,\eta}(t_c)}\,\tilde{w}_{j,m,l}(\frac{\xi}{\eta\,\g'(2^{-j})})\,\phi\left(\frac{\xi}{2^{m+j}}\right)\,
\phi^l\left(\frac{\eta}{2^{m+j}}\,\g'(\frac{1}{2^j})\right)\,,
\eeq
with $\tilde{w}_{j,m,l}\in C^{N-3}([2^{-k(\g)},\,2^{k(\g)}])$ and $\|\tilde{w}_{j,m,l}\|_{C^{N-3}([2^{-k(\g)},\,2^{k(\g)}])}\leq (C_{\g})^l$.

Let us now focus on the second summand in \eqref{mdecomp}. Proceeding as in the case of $\A_{j,m}$ we are writing
\beq\label{Brewrit}
\B_{j,m}^{k}(\xi,\eta)= \,2^{-\frac{m}{2}+k}\,e^{i\,\varphi_{\xi,\eta}(t_c)}\,b_{j,m,k,\eta}(t_c)\,\phi\left(\frac{\xi}{2^{m+j}}\right)\,
\phi\left(\frac{\eta}{2^{m+j}}\,\g'(\frac{1}{2^j})\right)\:,
\eeq
where for $u\in[2^{-k(\g)},\,2^{k(\g)}]$ we set
\beq\label{wdef}
b_{j,m,k,\eta}(u):=\int_{\R} \tilde{\va}(t)\,
e^{i\,\int_{0}^{2^{-\frac{m}{2}+k}\,t}\,\int_{0}^{s} \varphi''_{\xi,\eta}(r+u)\,dr\,ds}\:\r(2^{-\frac{m}{2}+k}t+u)\,dt\;.
\eeq
The idea here is to use the fact that the main domain of integration is further away from the critical point $t_c$.
With other words we will integrate by parts in \eqref{wdef}. For this, we define a differential operator $\tilde{L}$ for which
$$\tilde{L}(e^{i\,\int_{0}^{2^{-\frac{m}{2}+k}\,\cdot}\,\int_{0}^{s} \varphi''_{\xi,\eta}(r+t_c)\,dr\,ds})=e^{i\,\int_{0}^{2^{-\frac{m}{2}+k}\,\cdot}\,\int_{0}^{s} \varphi''_{\xi,\eta}(r+t_c)\,dr\,ds}\:.$$
Thus, we set
\beq\label{Ldif}
\tilde{L}:=\frac{1}{i 2^{-\frac{m}{2}+k} \int_{0}^{2^{-\frac{m}{2}+k}\,t}\, \varphi''_{\xi,\eta}(r+t_c)\,dr}\,\frac{\partial}{\partial t}
\eeq
with the transpose operator given by
\beq\label{Ldiftransp}
\tilde{L}^t h:= \frac{\partial}{\partial t}\left(\frac{i}{ 2^{-\frac{m}{2}+k} \int_{0}^{2^{-\frac{m}{2}+k}\,t}\, \varphi''_{\xi,\eta}(r+t_c)\,dr}\,h(\cdot)\right)\,
\eeq
Now, integrating by parts we have
$$b_{j,m,k,\eta}(t_c)=i\,\int_{\R}\frac{e^{i\,\int_{0}^{2^{-\frac{m}{2}+k}\,t}\,\int_{0}^{s} \varphi''_{\xi,\eta}(r+t_c)\,dr\,ds}}{2^{-\frac{m}{2}+k} \int_{0}^{2^{-\frac{m}{2}+k}\,t}\, \varphi''_{\xi,\eta}(r+t_c)\,dr}$$
$$\times (\tilde{\va}'(t)\, \r(2^{-\frac{m}{2}+k}t+t_c) +
\tilde{\va}(t)\,2^{-\frac{m}{2}+k}\,\r'(2^{-\frac{m}{2}+k}t+t_c))\,dt-$$
$$i\,\int_{\R}e^{i\,\int_{0}^{2^{-\frac{m}{2}+k}\,t}\,\int_{0}^{s} \varphi''_{\xi,\eta}(r+t_c)\,dr\,ds}\,\tilde{\va}(t)\,\r(2^{-\frac{m}{2}+k}t+t_c)\,\frac{\varphi''_{\xi,\eta}(2^{-\frac{m}{2}+k}t+t_c)}
{(\int_{0}^{2^{-\frac{m}{2}+k}\,t}\, \varphi''_{\xi,\eta}(r+t_c)\,dr)^2}\,dt\:.$$
Using \eqref{tcinv} and the fact that
\beq\label{expsrew}
\int_{0}^{2^{-\frac{m}{2}+k}\,t}\,\int_{0}^{s} \varphi''_{\xi,\eta}(r+t_c)\,dr\,ds=\int_{0}^{2^k t}\int_{0}^s \frac{\eta \g'(2^{-j})}{2^{m+j}}
\frac{\g''(\frac{r\,2^{-\frac{m}{2}}+t_c}{2^j})}{2^j\,\g'(2^{-j})}\,dr\,ds\
\eeq
one can check that
\beq\label{bjmk}
b_{j,m,k,\eta}(t_c)= \tilde{b}_{j,m,k} (\frac{\eta \g'(2^{-j})}{2^{m+j}},\,\frac{\xi}{\eta\,\g'(2^{-j})})\,,
\eeq
where $\tilde{b}_{j,m,k}(\cdot,\cdot):\, [\frac{1}{10},10]\times[2^{-k(\g)},\,2^{k(\g)}]\,\rightarrow\,\R$ is a $C^{N-3}$ function
with $\|\partial^{\a}\,\partial^{\b}\tilde{b}_{j,m,k}\|_{\infty}\lesssim_{\g} 2^{-2 k}$ if $\a+\b<N-3$ and $N$ large enough.

Putting together the information presented in \eqref{mdecomp}-\eqref{bjmk} we conclude that there exist the functions
$\{\zeta_k\}_{k=0}^{\frac{m}{2}+C_{\g}}$
with the properties
\beq\label{zeta}
\zeta_k:\,[\frac{1}{10},10]\times[2^{-k(\g)},\,2^{k(\g)}]\,\rightarrow\,\R\:\:\:\textrm{with}\:\:\:\;\:\:\|\zeta_k\|_{C^{N-3}}\lesssim_{\g} 2^{-k}
\eeq
such that
\beq\label{C}
m_{j,m,m}^{22}(\xi,\,\eta)= \sum_{k=0}^{\frac{m}{2}+C_{\g}} 2^{-\frac{m}{2}}\,e^{i\,\varphi_{\xi,\eta}(t_c)}\,
\zeta_k\left(\frac{\eta\,\g'(2^{-j})}{2^{m+j}},\,\frac{\frac{\xi}{2^{m+j}}}{\frac{\eta\,\g'(2^{-j})}{2^{m+j}}}\right)\,
\phi\left(\frac{\xi}{2^{m+j}}\right)\,\phi\left(\frac{\eta\,\g'(2^{-j})}{2^{m+j}}\right)\:,
\eeq
thus proving \eqref{centsymbol}.

\noindent Remark. Strictly speaking, each of the above functions $\zeta_k$ is in fact depending on $j,m$. However we have chosen not to write this explicit dependence, since the norm $\|\zeta_k\|_{C^N}$ is independent of the parameters $j$ and $m$.

\textbf{Observation.} Based on \eqref{zeta} and \eqref{C}, it will be enough to understand the behavior of the operator having
as multiplier the expression
 \beq\label{vjm}
v_{j,m}(\xi,\eta):= 2^{-\frac{m}{2}}\,e^{i\,\varphi_{\xi,\eta}(t_c)}\,
\zeta\left(\frac{\eta\,\g'(2^{-j})}{2^{m+j}},\,\frac{\frac{\xi}{2^{m+j}}}{\frac{\eta\,\g'(2^{-j})}{2^{m+j}}}\right)\,
\phi\left(\frac{\xi}{2^{m+j}}\right)\,\phi\left(\frac{\eta\,\g'(2^{-j})}{2^{m+j}}\right)
 \eeq
 where here we set $\zeta:=\zeta_0$.

 In Section 5, we will see that the function $\zeta$ appearing in the above expression, can be essentially replaced by the constant function $1$.

 Also notice that we have $\Psi_{\eta}(\xi)=-\varphi_{\xi,\eta}(t_c)$ where we recall that $\Psi_{\eta}$ is obtained from the relation $\Psi'_{\eta}(\cdot)=(\eta\,\g')^{-1}(\cdot)$.

$\newline$

From the above Observation, it remains to understand the main ``pieces"
\beq\label{pieces}
T_{j,m}(f,g)(x):= \int_{\R}\int_{\R}\hat{f}(\xi)\hat{g}(\eta)v_{j,m}(\xi,\eta) e^{i\xi x} e^{i\eta x} d\xi d\eta\:
 \eeq

The core of this paper will be how to show (and this constitutes the most difficult part of our result) that the operators $T_{j,m}$ obey
the condition
\beq\label{sumation}
\left\|\sum_{j\in\Z, m\geq 0}T_{j,m}(f,g)\right\|_1 \lesssim \left\|f\right\|_2\left\|g\right\|_2\:.
 \eeq
The plan for proving \eqref{sumation} will be described in the next section.

\section{\bf Elaborating on the key result}

From the above description our main theorem is reduced to the task of obtaining good bounds for each operator $T_{j,m}$. This is the point where we introduce one of the main novelties of the present paper: the ``scale type decay" of the $L^1-$norm of the operator $T_{j,m}$ relative to the level set of the multiplier's phase. More precisely, we claim that the following holds:

\begin{t1}\label{estimpiece}
There exists $\ep\in (0,1)$ such that
\beq\label{est}
\left\| T_{j,m}(f,g)\right\|_1\lesssim_{\g} 2^{-\ep m} \left\|f\right\|_2\left\|g\right\|_2\:.
\eeq
\end{t1}

This entire section will be dedicated for properly shaping the above theorem.

Firstly, notice that from the definition of $T_{j,m}$ relation \eqref{est} is equivalent with
 $$\left\| T_{j,m}(f,g)\right\|_1\lesssim_{\g} 2^{-\ep m} \left\|\hat{f}(\cdot)\,\phi\left(\frac{\cdot}{2^{m+j}}\right)\right\|_2\,
 \left\|\hat{g}(\cdot)\,\phi\left(\frac{\cdot\,\g'(2^{-j})}{2^{m+j}}\right)\right\|_2\:.$$
 Thus, once that we have proved \eqref{est}, one can follow the same reasonings as in the proof of Claim 2 - see \eqref{sumj}, and deduce that
\beq\label{sumjforT}
\left\| \sum_{j}T_{j,m}(f,g)\right\|_1\lesssim_{\g} 2^{-\ep m} \left\|f\right\|_2\left\|g\right\|_2\:.
\eeq

Let us now move to the preparatives for \eqref{est}.

From the Observations in Section 2, in what follows wlog we will assume that
$\lim_{t\rightarrow 0} \g'(t)=0$ and $\lim_{t\rightarrow\infty}\g'(t)=\infty$.

Next, while not necessary, we will choose to rescale our problem and thus make
the parallelism as well as the differences between our proof and Li's approach more transparent. Thus maintaining the notations from \cite{Li} we proceed as follows:

Using the scaling symmetry and depending on the values of $j$, we define the following operator:
\begin{itemize}

\item For $j>0$ (thus $2^{-j}\rightarrow 0$)

$$B_{j,m}(f(\cdot),g(\cdot))(x):=[\g'(2^{-j})]^{\frac{1}{2}}\,T_{j,m}\left( f(2^{m+j}\cdot), g(\frac{2^{m+j}}{[\g'(2^{-j})]}\cdot) \right)(\frac{x\,[\g'(2^{-j})]}{2^{m+j}})\:.$$
Remark that
$$B_{j,m}(f,g)(x)=2^{-\frac{m}{2}}\,[\g'(2^{-j})]^{\frac{1}{2}}\times$$
$$\int_{\R}\int_{\R}\hat{f}(\xi)\,\hat{g}(\eta)\,e^{i (\g'(2^{-j})\xi+\eta)\,x}\,e^{-i\,2^{m}\,2^j\,\Psi_{\frac{\eta}{\g'(2^{-j})}}(\xi)}\,
\zeta(\eta,\,\frac{\xi}{\eta})\,\phi(\xi)\,\phi(\eta)\,d\xi\,d\eta \:.$$
\noindent Remark. We used here the fact $\Psi_{a\,\eta}(a\,\xi)= a\,\Psi_{\eta}(\xi)$ for $a>0$.

At the \textit{heuristic} level, we should think at $B_{j,m}(f,g)(x)$ as given by the expression
\beq\label{heuristbjm}
\eeq
$$[\g'(2^{-j})]^{\frac{1}{2}}\,\int_{\R} \left(f*\check{\phi}\right)(\g'(2^{-j})x-2^{m}t)\,\left(g*\check{\phi}\right)(x+\frac{2^{m+j}}{\g'(2^{-j})}\g(\frac{t}{2^j}))\,\r(t)\,dt\:.$$
This heuristic will be helpful later in providing the intuition on \textit{how} to discretize our operator $B_{j,m}$.

\item For $j<0$ (thus $2^{-j}\rightarrow \infty$)

$$B_{j,m}(f(\cdot),g(\cdot))(x):=
[\g'(2^{-j})]^{-\frac{1}{2}}\,T_{j,m}\left( f(2^{m+j}\cdot), g(\frac{2^{m+j}}{[\g'(2^{-j})]}\cdot) \right)(\frac{x}{2^{m+j}})\:.$$
As before, notice that
$$B_{j,m}(f,g)(x)=2^{-\frac{m}{2}}\,[\g'(2^{-j})]^{-\frac{1}{2}}\times$$
$$\int_{\R}\int_{\R}\hat{f}(\xi)\,\hat{g}(\eta)\,e^{i (\xi+\frac{\eta}{[\g'(2^{-j})]})\,x}\,e^{-i\,2^{m}\,2^j\,\Psi_{\frac{\eta}{\g'(2^{-j})}}(\xi)}\,
\zeta(\eta,\,\frac{\xi}{\eta})\,\phi(\xi)\,\phi(\eta)\,d\xi\,d\eta \:.$$
Again, at the heuristic level, we have that $B_{j,m}$ encapsulates the behavior of
\beq\label{heuristbjm1}
\eeq
$$[\g'(2^{-j})]^{-\frac{1}{2}}\,\int_{\R} \left(f*\check{\phi}\right)(x-2^{m}t)\,\left(g*\check{\phi}\right)(\frac{x}{\g'(2^{-j})}+\frac{2^{m+j}}{\g'(2^{-j})}\g(\frac{t}{2^j}))\,\r(t)\,dt\:.$$
\end{itemize}

With these definitions we see that the statement of Theorem 2 is equivalent with

\textbf{Theorem 2 (reformulated).}
\textit{There exists $\ep\in (0,1)$ such that}
\beq\label{estref}
\left\| B_{j,m}(f,g)\right\|_1\lesssim_{\g} 2^{-\ep m} \left\|f\right\|_2\left\|g\right\|_2\:.
\eeq

\section{\bf Discretization of our operator $B_{j,m}$}

We start by saying that we will focus on the case $j>0$ as the reasonings for the other case ($j<0$) are of similar nature.
Thus in what follows one should think at $2^{-j}$ as being a very small quantity.

We first turn our attention towards the Fourier side of our operator and remodel the main term of its phase:
\beq\label{Fourierbjm}
\eeq
$$B_{j,m}(f,g)(x)=2^{-\frac{m}{2}}\,[\g'(2^{-j})]^{\frac{1}{2}}$$
$$\int_{\R}\int_{\R}\hat{f}(\xi)\,\hat{g}(\eta)\,e^{i (\g'(2^{-j})\xi+\eta)\,x}\,e^{-i\,2^{m}\,2^j\,\Psi_{\frac{\eta}{\g'(2^{-j})}}(\xi)}\,
\zeta(\eta,\,\frac{\xi}{\eta})\,\phi(\xi)\,\phi(\eta)\,d\xi\,d\eta \:.$$

Now based on \eqref{asymptotic0gaminv}, for $x\in J$,  we have that
\beq\label{rel1}
2^{j}\,(\g')^{-1}(x\,\g'(2^{-j}))= r(x)\,+\,r_j(x)\:.
\eeq

Using the above relation, we notice that up to a factor depending possibly only of $j$ and $\eta$ that can be absorbed in $\hat{g}$, for $\xi,\:\eta\in I$ with $\frac{\xi}{\eta}\in J$ we have
\beq\label{Psi}
\begin{array}{rl}
&2^j\,\Psi_{\frac{\eta}{\g'(2^{-j})}}(\xi)=\int_{\eta}^{\xi}2^j\,(\g')^{-1}(\frac{u}{\eta}\,\g'(2^{-j}))\,du\\
 &=\int_{\eta}^{\xi} r(\frac{u}{\eta})\, du\,+\,\int_{\eta}^{\xi} r_j(\frac{u}{\eta}):=
 \eta\,R(\frac{\xi}{\eta})\,+\,\eta\,R_j(\frac{\xi}{\eta})\,.
\end{array}
\eeq
Now since the $R_j$'s are smooth and tend uniformly to zero we do expect that their contribution is negligible in the behavior of $B_{j,m}$; in what follows, for notational simplicity we will only write/keep the dominant oscillatory factor, namely the exponential given by $R$. Alternatively, the reader can check that each of our reasonings below carry over smoothly in the presence of the error term in \eqref{Psi} (all the estimates for $R$ remain true for $R+R_j$).

 Thus, we have obtained
\beq\label{Bjmintdec}
\eeq
$$B_{j,m}(f,g)(x)=
 2^{-\frac{m}{2}}[\g'(2^{-j})]^{\frac{1}{2}} $$
 $$\int_{\R}\int_{\R}\hat{f}(\xi)\,\hat{g}(\eta)\,e^{i (\g'(2^{-j})\xi+\eta)\,x}\,
e^{-i\,2^m\,\eta\,R(\frac{\xi}{\eta})}\,\zeta(\eta,\,\frac{\xi}{\eta})\,\phi(\xi)\,\phi(\eta)\,d\xi\,d\eta \:.$$
$\newline$

\noindent\textbf{5.1. Discretized model - first step: flattening the function $\z$.}
$\newline$

The message of this subsection is that, due to the smoothness of $\z$, will be enough to treat the expression of $B_{j,m}$
 in which we replace the function $\z$ with the constant function $1$. To see this, we first notice that wlog
 we may suppose that $\zeta$ is compactly supported in $[0,2\pi]\times[0,2\pi]$. Looking now at it as a periodic function and
 applying standard Fourier analysis we have that
\beq\label{fseries}
\zeta(\eta,\frac{\xi}{\eta})=\sum_{n_1,\,n_2\in\Z} c_{k_1,\,k_2}\,e^{i\,\eta\,k_1}\, e^{i\,\frac{\xi}{\eta}\,k_2}\,,
\eeq
where from the hypothesis on $\z$, namely that $\|\zeta\|_{C^N}\lesssim_{\g} 1$ we also have
\beq\label{coefdecay}
|c_{k_1,\,k_2}|\lesssim_{\g} \frac{1}{(|k_1|+|k_2|+1)^{4}}\:.
\eeq
Thus, we deduce that
\beq\label{bdec}
B_{j,m}(f,g)=\sum_{k_1,\,k_2\in\Z} c_{k_1,\,k_2}\, B_{j,m}^{k_1,k_2}(f,g)
\eeq
with
\beq\label{bjmk1k2}
\eeq
$$B_{j,m}^{k_1,k_2}(f,g):=2^{-\frac{m}{2}}\,[\g'(2^{-j})]^{\frac{1}{2}}$$
$$\int_{\R}\int_{\R}\hat{f}(\xi)\,\hat{g}(\eta)\,e^{i (\g'(2^{-j})\xi+\eta)\,x}\,e^{-i\,2^m\,\eta\,R(\frac{\xi}{\eta})}\,
\,e^{i\,\eta\,k_1}\, e^{i\,\frac{\xi}{\eta}\,k_2}\,\phi(\xi)\,\phi(\eta)\,d\xi\,d\eta\:. $$
As a direct consequence of the methods that will be exposed at the Step 3 in our discretization process as well as of those in the next section
we have that the following extension of Theorem 2 holds:

\begin{t1}\label{th2extended}
There exist $\ep\in (0,1)$ and $C^1_\g>0$ small enough (any $C^1_\g$ smaller than $10^{-4}\,c_{\g}$ with $c_{\g}$ defined in \eqref{fstterm0} works) such that for any $|k_2|\leq C^1_{\g}\,2^{m}$ the following holds:
\beq\label{bjmk1k2estsmall}
\left\| B_{j,m}^{k_1,k_2}(f,g)\right\|_1\lesssim_{\g} 2^{-\ep m} \left\|f\right\|_2\left\|g\right\|_2\:.
\eeq
\end{t1}

\textbf{Observations.} 1) Notice that the uniform estimates in the parameter $k_1$ are trivial
since the factor $e^{i\,\eta\,k_1}$ can be absorbed in the function $\hat{g}(\eta)$, without changing the nature of our bounds.

2) As mentioned before, Theorem 3 is just a direct consequence of the proof of our Theorem 2. As a result, we will not provide here any details
but only stress that all the arguments following this first stage of discretization are easily adapted to the context in which the multiplier of our operator contains the extra-factor $ e^{i\,\frac{\xi}{\eta}\,k_2}$. Indeed all that one needs to is that the $\xi$-\textit{derivative} of the main term $2^m\,\eta\,R(\frac{\xi}{\eta})$ of the phase remains dominant when adding the term $\frac{\xi}{\eta}\,k_2$ (this is required when proving the analogue of \eqref{secondmomemtcond}). This is precisely where the condition $|k_2|\leq C^1_{\g}\,2^{m}$ with $C^1_{\g}$ small intervenes. This condition assures that the key estimates \eqref{secondmomemtcond}, \eqref{critreg}, \eqref{critreg} and \eqref{hormandphase} remain valid. We leave further details to the reader.

For the large values of $|k_2|$ we have the following simpler to check result
$\newline$

\noindent\textbf{Proposition}. \textit{If $C_{\g}>0$ is any positive constant and $|k_2|>C_{\g}\,2^{m}$ then}
\beq\label{bjmk1k2estlarge}
\|B_{j,m}^{k_1,k_2}(f,g)\|_{1}\lesssim_{\g} |k_2|^{\frac{3}{2}}\,\|f\|_2\,\|g\|_2\;.
\eeq
\noindent Remark. Here one could obtain better estimates; however this statement is enough for our proof.

\begin{proof}
Here is a sketch of the proof.

Following 1) in the Observation above we may assume wlog that $k_1=0$.

Set \beq\label{frame}
\varphi_{m,l,p}(x)= 2^{-\frac{m}{2}}\,\varphi(2^{-m}x-l)\,e^{i\,2^{-m}\,p\,x}
\eeq
with $\varphi$ smooth such that $\textrm{supp}\,\hat{\varphi}\subseteq [0,2]$.

Also, for notational simplicity, let us define in what follows $\v_{0,l}:=\v_{0,l,0}\:.$

Using now \eqref{bjmk1k2}, we decompose
\beq\label{decfg}
\eeq
$$f=\sum_{l\in \Z} <f,\,\varphi_{0,l}>\,\v_{0,l}\,,$$
$$g=\sum_{l'\in \Z} <g,\,\varphi_{0,l'}>\,\v_{0,l'}\,.$$

\noindent Remark. Here we are less ambitious with our discretization (coarser frequency scale) since we only intend to obtain some relaxed bounds on $B_{j,m}$.

Then
$$B_{j,m}^{k_2}(f,g)(x)=2^{-\frac{m}{2}}\,[\g'(2^{-j})]^{\frac{1}{2}}\,\sum_{l,l'\in\Z} <f,\,\varphi_{0,l}> <g,\,\varphi_{0,l'}> A_{l,l'}^{k_2}(x)$$
with
\beq\label{All}
 A_{l,l'}^{k_2}(x):=\int_{\R}\int_{\R}\,e^{i\,[-\xi l-\eta l'+(\g'(2^{-j})\xi+\eta)\,x-2^m\,\eta\,R(\frac{\xi}{\eta})+\frac{\xi}{\eta}\,k_2]}\,
\hat{\v}(\xi)\,\hat{\v}(\eta)\,\phi(\xi)\,\phi(\eta)\,d\xi\,d\eta \:.
\eeq
Of course the central element will be to take advantage of the oscillation of the phase. As always the main contributions arise from the regions
near the stationary points.

Choosing now $\nu\in C_0^{\infty}(\R)$ with $\textrm{supp}\,\nu\subset[-10,10]$ and
$$\sum_{r\in\Z} \nu(x-r)=1\;\:\:\:\;\:\:\:\:\:\forall\:x\in\R\,,$$
we write
\beq\label{comppart}
\eeq
$$B_{j,m}^{k_2}(f,g)(x)=2^{-\frac{m}{2}}\,[\g'(2^{-j})]^{\frac{1}{2}}\,$$
$$\sum_{r,r'\in\Z}\sum_{l,l'\in\Z} <f,\,\varphi_{0,l}> <g,\,\varphi_{0,l'}> A_{l,l'}^{k_2}(x)\,\nu(k_2^{-1}\,\g'(2^{-j})\,x-r)\,\nu(k_2^{-1}\,x-r')\;.$$
Applying now the principle of (non)stationary phase we notice that
\begin{itemize}
\item{If $|r\,k_2-l|\gtrsim_{\g}|k_2|$ then
$$|A_{l,l'}^{k_2}(x)\,\nu(k_2^{-1}\,\g'(2^{-j})\,x-r)|\lesssim_{\g}\frac{1}{|r\,k_2-l|^{10}}\:.$$}
\item{If $|r'\,k_2-l'|\gtrsim_{\g}|k_2|$ then
$$|A_{l,l'}^{k_2}(x)\,\nu(k_2^{-1}\,x-r')|\lesssim_{\g}\frac{1}{|r'\,k_2-l'|^{10}}\:.$$}
\item{For any $l,l'\in\Z$ one has
$$|A_{l,l'}^{k_2}(x)|\lesssim_{\g} \min\{(2^{m})^{-\frac{1}{2}},\,(k_2)^{-\frac{1}{2}}\}\:.$$}
\end{itemize}
Thus the main term in \eqref{comppart} is represented by the $l,l'$ summation in the range $|r\,k_2-l|\lesssim_{\g}|k_2|$ and
$|r'\,k_2-l'|\lesssim_{\g}|k_2|$ respectively.

Finally, \eqref{bjmk1k2estlarge} follows now by applying a Cauchy-Schwarz argument to the main term.
\end{proof}

Combining now  \eqref{bjmk1k2estsmall} and \eqref{bjmk1k2estlarge} with \eqref{coefdecay} and \eqref{bdec} we conclude that we can reduce
the study of $B_{j,m}$ to that of $B_{j,m}^{0,0}$. For notational simplicity we will \textbf{re-denote} $B_{j,m}=B_{j,m}^{0,0}$.

$\newline$
\textbf{Remark.} As claimed previously, we have now reduced our problem to understanding the initial $B_{j,m}$ operator in the particular case $\zeta=1$:
\beq\label{Bjmreduction}
\eeq
$$B_{j,m}(f,g)(x)=
 2^{-\frac{m}{2}}[\g'(2^{-j})]^{\frac{1}{2}} $$
 $$\int_{\R}\int_{\R}\hat{f}(\xi)\,\hat{g}(\eta)\,e^{i (\g'(2^{-j})\xi+\eta)\,x}\,
e^{-i\,2^m\,\eta\,R(\frac{\xi}{\eta})}\,\phi(\xi)\,\phi(\eta)\,d\xi\,d\eta \:.$$
$\newline$

\noindent\textbf{5.2. Discretized model - second step: gaining intuition about how to split our operator.}
$\newline$

Before we should proceed with our discretization it is worth mentioning several aspects that will help our understanding
of the operator $B_{j,m}$.

Following \eqref{heuristbjm} in which we make use of property \eqref{asymptotic0}, we have the following
$\newline$
\textbf{Main Heuristic:}
\beq\label{heuristbjm2}
\eeq
$$B_{j,m}(f,g)(x)\approx$$
$$[\g'(2^{-j})]^{\frac{1}{2}}\,\int_{\R} \left(f*\check{\phi}\right)(\g'(2^{-j})x-2^{m}t)\,
\left(g*\check{\phi}\right)(x+2^{m}\,Q(t))\,\r(t)\,dt\:.$$

Thus, based on this remark, we make the following

$\newline$
\textbf{Observations:}
\begin{itemize}
\item when $t$ varies inside $\textrm{supp}\,\r$, the space variable of $g$ covers an interval of length $2^m$ (heuristically this means averaging the value of $g$ along intervals of length $2^m$).
\item for the space variable $(\g'(2^{-j})x-2^{m}t)|_{t\in\textrm{supp}\,\r}$ corresponding to $f$, to cover an interval of length $2^m$, the $x$-variable is allowed to move in an interval of length $\frac{2^{m}}{\g'(2^{-j})}$.
\end{itemize}
These observations have the following
$\newline$

\noindent\textbf{Implications:}
\begin{itemize}
\item{the first observation suggests for $g$ a Gabor-decomposition of the form
\beq\label{gabor}
g(x)=\sum_{l,p} \left\langle g, \varphi_{m,l,p} \right\rangle\:\varphi_{m,l,p}\;,
\eeq
where here we have used the same notations as in \eqref{frame}.}
\item{the second observation has the following \textbf{consequence}:
\beq\label{localization}
\eeq
\textit{It is enough to prove \eqref{estref} for $B_{j,m}(f,g)$ and $g$ restricted to the interval $[0,\frac{2^{m}}{\g'(2^{-j})}]$.}
$\newline$
Indeed, at the heuristic level, our $B_{j,m}$ is a \textbf{local operator}. To be more precise we have the following:
 set $I_{m,j}^{u}:=[\frac{2^m\,u}{|\g'(2^{-j})|},\,\frac{2^m\,(u+1)}{|\g'(2^{-j})|}]$ and $I_{m}^{s}:=[2^m\,s,\,2^m\,(s+1)]$ where $u,\,s\in\Z$;
 then, for $v\in\Z$ and $|u-v|,\,|v-s|>>C_{\g}$ we have
 \beq\label{heuristiclocaliz}
 \|B_{j,m}(f\,\chi_{I_{m}^s}, g\,\chi_{I_{m,j}^u})\|_{L^1(I_{m,j}^v)}\lesssim_{\g} \frac{2^{-m}}{|u-v|^3\,|v-s|^3}\,
\|f\|_{L^2(I_{m}^s)}\,\|g\|_{L^2(I_{m,j}^u)}\:.
 \eeq
Here the notation $\chi_{I_{m}^s}$ designates a smooth version of the characteristic function of the interval $I_{m}^s$ (similarly for $\chi_{I_{m,j}^u}$). The justification of \eqref{heuristiclocaliz} will become apparent at the end of Stage 1, third step below, by just inspecting  \eqref{Bjm} and \eqref{Qmpdef}.
}
\end{itemize}

Now, based on the first Implication above we will write
\beq\label{Bjmreductionn0}
\eeq
$$B_{j,m}(f,g)(x)=\sum_{l\in\Z}\sum_{p=2^m}^{2^{m+1}}
 2^{-\frac{m}{2}}[\g'(2^{-j})]^{\frac{1}{2}} \left\langle g, \varphi_{m,l,p} \right\rangle\,$$
 $$\int_{\R}\int_{\R}\hat{f}(\xi)\,\hat{\varphi}_{m,l,p}(\eta)\,e^{i (\g'(2^{-j})\xi+\eta)\,x}\,
e^{-i\,2^m\,\eta\,R(\frac{\xi}{\eta})}\,\phi(\xi)\,\phi(\eta)\,d\xi\,d\eta \:.$$

\noindent Remark. The summation restriction in the $p$-parameter is a result of the compact support of $\phi$. Moreover, as a consequence of \eqref{localization} and \eqref{heuristiclocaliz} we will later reduce the study of our operator to a modified
\eqref{Bjmreductionn0} where this time the $l$ parameter only ranges between $0$ and $[\g'(2^{-j})]^{-1}$.

$\newline$
\noindent\textbf{5.3. Discretized model - third step: separation of variables.}
$\newline$

The key idea in this subsection is to decompose the multiplier of $B_{j,m}$ as a superposition of tensor products.
$\newline$

\noindent\textbf{Stage 1. Phase discretization.}
$\newline$

We continue our program, by further dicretizing our operator, namely by freezing the variable $\eta$ of the phase represented in \eqref{Bjmreduction}.

Indeed, we first notice that for $\xi\in\textrm{supp}\:\phi$  and $\eta\in \textrm{supp}\,\hat{\varphi}_{m,l,p}$ we have

$$2^{m}\,\left|\eta\,R(\frac{\xi}{\eta})-\frac{p}{2^m}\,R(\frac{\xi}{\frac{p}{2^m}}) \right|\lesssim_{\g} 1\:.$$

Canonical reasonings by now, see below, will allow us to further reduce the analysis of $B_{j,m}$ to the study of
\beq\label{Bjm}
\begin{array}{cl}
\B_{j,m}(f,g)(x)=\\
\sum_{l\in\Z}\sum_{p=2^m}^{2^{m+1}}
 2^{-\frac{m}{2}}[\g'(2^{-j})]^{\frac{1}{2}} \left\langle g, \varphi_{m,l,p} \right\rangle\,\Q_{m,p}f(\g'(2^{-j})\,x)\,\varphi_{m,l,p}(x)\:,
\end{array}
\eeq
where
\beq\label{Qmpdef}
 \Q_{m,p}f(x):=\int_{\R}\hat{f}(\xi)\,\phi(\xi)\,e^{-i\,p\,R(\frac{\xi}{\frac{p}{2^m}})}\,e^{i\,\xi\,x}\,d\xi\:.
\eeq
Indeed, to see that this is enough, we first notice that for $\eta, \eta_0\in \textrm{supp}\:\varphi_{m,l,p}$  (hence $|2^{m}(\eta-\eta_0)|\lesssim 1$) and $h(\xi,\eta):=\frac{\eta\,R(\frac{\xi}{\eta})-\eta_0\,R(\frac{\xi}{\eta_0})}{\eta-\eta_0}$ we have $h\in C^{N-1}$ and thus for any $k\in\N$ we deduce $h^k(\xi,\eta)=\sum_{n_1,n_2\in\Z} d_{k,n_1,n_2}\,e^{i\xi\,n_1}\,e^{i\eta\,n_2}\,,$
with $$|d_{k,n_1,n_2}|\lesssim_{r} \frac{k^{r}}{|n_1|^r+|n_2|^r}\:\|h\|_{C^N}^{k}\:\:\:\:\:\:\:\:\:\:\:\:\:\:\:\:\:\forall\:0\leq r\leq N-1.$$

Then, using the Taylor expansion
$$e^{i\,2^m (\eta-\eta_0)\,h(\xi,\eta)}=\sum_{k\in\N}\frac{[2^m (\eta-\eta_0)]^k}{k!}\,i^k\,\sum_{n_1,n_2\in\Z} d_{k,n_1,n_2}\,e^{i\xi\,n_1}\,e^{i\eta\,n_2}\,$$
and the decay of the coefficients $d_{k,n_1,n_2}$ we have that indeed $\B_{j,m}$ can be reduced to a discretized version of type \eqref{Bjm}.

At this point setting $B_{j,m}$ as in \eqref{Bjm}, and based on \eqref{Qmpdef}, one can rigorously justify formulation \eqref{heuristiclocaliz}.
Thus, from now on, we will refer to $B_{j,m}$ as given by

\beq\label{Bjm1}
\begin{array}{cl}
B_{j,m}(f,g)(x)=\\
\sum_{l=0}^{[\g'(2^{-j})]^{-1}}\sum_{p=2^m}^{2^{m+1}}
 2^{-\frac{m}{2}}[\g'(2^{-j})]^{\frac{1}{2}} \left\langle g, \varphi_{m,l,p} \right\rangle\,\Q_{m,p}f(\g'(2^{-j})\,x)\,\varphi_{m,l,p}(x)\:.
\end{array}
\eeq

$\newline$

\noindent\textbf{Stage 2. Local flattening of the operator $\Q_{m,p}f$ .}
$\newline$

Remark now that we have
\beq\label{qmp}
\left|\Q_{m,p}f(x)-\Q_{m,p}f(x')\right|\lesssim |x-x'|\,\|f\|_2\:.
\eeq
This implies that $\Q_{m,p}f$ ``sees'' the exterior environment in unit steps. As a consequence,
will become natural to apply two types of discretizations depending on the relative sizes of the ``moral" support of $\Q_{m,p}f(\g'(2^{-j})\,\cdot)$ and of the support of $\varphi_{m,l,p}(\cdot)$.

\textbf{Claim.} For proving \eqref{estref}, it is enough to show the analogous discretized version
\beq\label{estrefdiscret}
\left\| \tilde{B}_{j,m}(f,g)\right\|_1\lesssim_{\g} 2^{-\ep m} \left\|f\right\|_2\left\|g\right\|_2\:.
\eeq
where the original $B_{j,m}$ is replaced by the model operator $\tilde{B}_{j,m}$ given by:
\begin{itemize}
\item \textbf{Case 1.} If $\g'(2^{-j})>2^{-m}$ then
$$\tilde{B}_{j,m}(f,g)(x)=2^{-m}\sum_{l=0}^{[\g'(2^{-j})]^{-1}}\sum_{p=2^m}^{2^{m+1}}\sum_{k=0}^{[2^{m}\g'(2^{-j})]}
 \left\langle g, \varphi_{m,l,p} \right\rangle\,$$ $$\Q_{m,p}f(l\,2^{m}\g'(2^{-j})+k)\,
 \tilde{\varphi}_{\log_2 \frac{1}{\g'(2^{-j})},l2^{m}\g'(2^{-j})+k,\frac{p}{2^{m}\,\g'(2^{-j})}}(x)\:.$$
\item \textbf{Case 2.} If $\g'(2^{-j})\leq 2^{-m}$ then
$$\tilde{B}_{j,m}(f,g)(x)=$$
 $$2^{-\frac{m}{2}}[\g'(2^{-j})]^{\frac{1}{2}}\sum_{l=0}^{[\g'(2^{-j})]^{-1}}\sum_{p=2^m}^{2^{m+1}}
 \left\langle g, \varphi_{m,l,p} \right\rangle\,\Q_{m,p}f(l\,2^{m}\g'(2^{-j}))\,\tilde{\varphi}_{m,l,p}(x)\:.$$
\end{itemize}

\noindent Remark. The function $\tilde{\varphi}$ above is only
assumed to be a Schwartz function. As before we set $\tilde{\varphi}_{m,l,p}(x)= 2^{-\frac{m}{2}}\,\tilde{\varphi}(2^{-m}x-l)\,e^{i\,2^{-m}\,p\,x}$.
$\newline$

\noindent\textbf{Proof of our Claim.}
$\newline$

In what follows we will only focuss on the second case, i.e. $\g'(2^{-j})\leq 2^{-m}$. The first case has a similar treatment and is only sketched
in the Observations (fourth one) ending this section.

Choose now $\u\in C_0^{\infty}(\R)$ with $\textrm{supp}\,\u\subset [-10,10]$ such that
\beq\label{partunity}
1=\sum_{s\in\Z} \u(x-s)\:\:\;\:\:\:\;\:\:\:\:\:\:\:\:\:\:\:\forall\:x\in\R\:.
\eeq

Using \eqref{partunity}, for each $l\in \{0,\ldots,[\g'(2^{-j})]^{-1}\}$ we write
$$\Q_{m,p}f(\g'(2^{-j})\,x)=\sum_{s\in\Z} \Q_{m,p}f(\g'(2^{-j})\,x)\,\u(2^{-m}\,x-l-s)=$$
$$\sum_{s\in\Z} (\int_{\R}\hat{f}(\xi)\,\phi(\xi)\,e^{-i\,p\,R(\frac{\xi}{\frac{p}{2^m}})}\,
e^{i\,\xi\,\g'(2^{-j})\,2^m\,(l+s)}\,e^{i\,\xi\,\g'(2^{-j})\,2^m\,(2^{-m}\,x-l-s)}\,d\xi)\,\u(2^{-m}\,x-l-s)=$$

$$=\sum_{s\in\Z}\sum_{k=0}^{\infty}\frac{(i\,\g'(2^{-j})\,2^m)^k}{k!}\,
\left(\int_{\R}\hat{f}(\xi)\,\xi^k\,\phi(\xi)\,e^{-i\,p\,R(\frac{\xi}{\frac{p}{2^m}})}\,
e^{i\,\xi\,\g'(2^{-j})\,2^m\,(l+s)}\,d\xi\right)$$
$$\times\,(2^{-m}\,x-l-s)^k\,\u(2^{-m}\,x-l-s)\;.$$

Now, for $k\in\N$, we set $\phi^{k}(x)=x^k\,\phi(x)$ and $\u^{k}(x):=x^k\,\u(x)$ respectively. As before, we let $\u^{k}_{m,l}(x):=\u^{k}(2^{-m}\,x-l)$.

Further, define
$$\Q_{m,p}^{k}f(x):=\int_{\R}\hat{f}(\xi)\,\phi^{k}(\xi)\,e^{-i\,p\,R(\frac{\xi}{\frac{p}{2^m}})}\,e^{i\,\xi\,x}\,d\xi\:.$$

With the previous notations, define now
\beq\label{discretbjm0}
\eeq
$$\tilde{B}^{s,k}_{j,m}(f,g)(x):= $$
$$2^{-\frac{m}{2}}\,[\g'(2^{-j})]^{\frac{1}{2}}\,\sum_{l=0}^{[\g'(2^{-j})]^{-1}}\,\sum_{p=2^m}^{2^{m+1}}
 \left\langle g, \varphi_{m,l,p} \right\rangle\,\Q^{k}_{m,p}f((l+s)\,2^{m}\g'(2^{-j}))\,\varphi_{m,l,p}(x)\,\u^{k}_{m,l+s}(x)\:,$$
and
$$\tilde{B}_{j,m}^{s}(f,g)(x)=\sum_{k=0}^{\infty} \frac{(i\,\g'(2^{-j})\,2^m)^k}{k!}\,\tilde{B}^{s,k}_{j,m}(f,g)(x)\:.$$

Based on the above considerations we deduce that

\beq\label{discretbjm}
\B_{j,m}(f,g)(x)=\sum_{s\in\Z} \tilde{B}_{j,m}^{s}(f,g)(x)\:.
\eeq

\textbf{Observations.}

1) From \eqref{discretbjm0}, the assumption $|\g'(2^{-j})|\,2^{m}\leq 1$ and the key fact that $\phi,\,\u\in C_0^{\infty}(\R)$ we deduce that the summation in $k$ is under control - absolutely summable series.

2) In the next section, using a duality argument, we will provide bounds on the $L^1$ norm of $\tilde{B}_{j,m}^{0}(f,g)$ for $f,\,g\in L^2$. A simple inspection of that proof gives us the following conclusion: the operator norm $\|\cdot\|_{L^2\times L^2\rightarrow L^1}$ for the corresponding (bilinear) operator $\tilde{B}^{s,k}_{j,m}(\cdot,\cdot)$ will only depend on the constants involving the properties of the curve $\g$ and on the following expressions:
\beq\label{normdep}
\eeq
$$\sup_{\|f\|_{L^2}\leq 1}\|\hat{f}(\xi)\,\phi^{k}(\xi)\|_{2}\leq 100^k\,,$$
$$\sup_{{{a\in\N,\,b\in \N^{*}}\atop{a,b\lesssim [\g'(2^{-j})]^{-1}}}\atop{\|h\|_{\infty}\leq1}}\left\{\left(\frac{2^{-m}\,\g'(2^{-j})}{b}\right)^{\frac{1}{2}}\:
\sum_{l=a}^{a+b} \|h\,\varphi_{m,l,0}\:\u^{k}_{m,l+s}\|^2_{2}\right\}^{\frac{1}{2}}\lesssim_N2^{-\frac{m}{4}}\,\frac{100^k}{1+s^{2N}}\:.$$

3) As a consequence, one deduces that
\beq\label{estrefdiscretk}
\left\| \tilde{B}_{j,m}^{s}(f,g)\right\|_1\lesssim_{\g, N} \frac{2^{-\ep m}}{1+s^{2N}} \left\|f\right\|_2\left\|g\right\|_2\:.
\eeq
Thus, from  \eqref{discretbjm} and \eqref{estrefdiscretk} our claim is proved for the Case 2.

4) For the first case, one follows the same lines - partition of unity and Taylor expansions. With the previous notations,
one decomposes $B_{j,m}$ as in \eqref{discretbjm}, where this time
\beq\label{discretcase1}
\tilde{B}_{j,m}^{s}(f,g)=\sum_{r\in\N} \frac{i^r}{r!}\,\tilde{B}_{j,m}^{s,r}(f,g)\,,
\eeq
with
$$\tilde{B}_{j,m}^{s,r}(f,g)(x)=2^{-m}\sum_{l=0}^{[\g'(2^{-j})]^{-1}}\sum_{p=2^m}^{2^{m+1}}\sum_{k=0}^{[2^{m}\g'(2^{-j})]}
 \left\langle g, \varphi_{m,l,p} \right\rangle\,\Q_{m,p}^{r}f(l\,2^{m}\g'(2^{-j})+k+s)$$ $$
 \tilde{\varphi}_{\log_2 \frac{1}{\g'(2^{-j})},l2^{m}\g'(2^{-j})+k,\frac{p}{2^{m}\,\g'(2^{-j})}}(x)\,\:\:\u^{r}_{\log_2 \frac{1}{\g'(2^{-j})},l2^{m}\g'(2^{-j})+k+s}(x)\:.$$
With these facts, following line by line the steps presented in the proof of Proposition 1, for the regime $2^m\,\g'(2^{-j})$ ``small"
one finds that
\beq\label{estrefdiscretk1}
\left\| \tilde{B}_{j,m}^{s}(f,g)\right\|_1\lesssim_{\g, N} \frac{2^{-\frac{m}{4}}}{1+s^{2N}}\,(2^m\,\g'(2^{-j}))^{\frac{1}{2}}\, \left\|f\right\|_2\left\|g\right\|_2\:.
\eeq
This ends the discretization of our operator.

\section{\bf The proof of the main theorem}

In this section we will prove the boundedness of the bilinear Hilbert transform along the curve $\Gamma$.
For notational simplicity we will drop the $\tilde{}$ symbol from both the definition of $\tilde{B}_{j,m}$ and $\tilde{\varphi}$.

\noindent Remark. In what follows, we will only use that the functions $\varphi,\,\tilde{\varphi}$ are Schwartz; the relative
localization of these functions is irrelevant and hence identifying $\varphi$ with $\tilde{\varphi}$ is safe.

Define
\beq\label{duality}
\l_{j,m}(f,g,h):=\int_{\R}B_{j,m}(f,g)(x)\, \overline{h(x)}\, dx=\int_{\R}f(x)\,\overline{D_{j,m}(g,h)(x)}\, dx\:.
\eeq
Then, Theorem 2 is a direct consequence of the following

\begin{t1} There exists $\ep\in (0,1)$ such that
$$|\l_{j,m}(f,g,h)|\lesssim_{\g} 2^{-\ep m}\|f\|_2\|g\|_2\|h\|_{\infty}\:.$$
\end{t1}
Further, the above theorem follows from

\begin{p1}\label{tt*}
For $\g'(2^{-j})>2^{-m}$ we have that
$$\left\| B_{j,m}(f,g)\right\|_1\lesssim_{\g} 2^{-\frac{m}{16}}\left\|f\right\|_2\left\|g\right\|_2\:.$$
\end{p1}

\begin{p1}
For $\g'(2^{-j})\leq 2^{-m}$ we have that
$$\|D_{j,m}(g,h)\|_2\lesssim_{\g} 2^{-\frac{m}{4}}\|g\|_2\|h\|_{\infty}\:.$$
\end{p1}

$\newline$
\textbf{Proof of Proposition 2.}
$\newline$

Making use of the discretization of $B_{j,m}$ for $\g'(2^{-j})\leq 2^{-m}$ described in Section 4, and using the dual expression
\eqref{duality}, we have
$$D_{j,m}(g,h)=2^{-\frac{m}{2}} [\g'(2^{-j})]^{\frac{1}{2}} \sum_{l=0}^{[\g'(2^{-j})]^{-1}}\sum_{p=2^m}^{2^{m+1}}
 \overline{\left\langle g, \varphi_{m,l,p} \right\rangle}\,
 \left\langle h, \varphi_{m,l,p} \right\rangle\,\psi_{m,l,p}(x)\:,$$
 where $\overline{\hat{\psi}_{m,l,p}(\xi)}:=\,e^{-i\,p\,R(\frac{\xi}{\frac{p}{2^m}})}\,e^{i\,l\,2^{m}\,\g'(2^{-j})\,\xi}\,\phi(\xi)\:.$

Now, the first natural step is to understand the interaction
 \beq\label{interact}
\left\langle\psi_{m,l,p},\psi_{m,l',p'} \right\rangle=
\int_{\R} e^{i\,[p\,R(\frac{\xi}{\frac{p}{2^m}})-p'\,R(\frac{\xi}{\frac{p'}{2^m}})]}\,
e^{-\,i\,\g'(2^{-j})(l-l')\,2^{m}\xi}\,\phi(\xi)\,\overline{\phi(\xi)}\,d\xi\:.
\eeq
 Applying the change of variable $\xi\rightarrow \frac{t}{2^m}$ and using the stationary phase principle  with $\Phi(t)=\Phi_{p,p'}(t):=p\,R(\frac{t}{p})-p'\,R(\frac{t}{p'})$ (and $\Psi=\Psi_{p,p'}$ the corresponding dual function) we deduce that
\beq\label{interactestim}
\begin{array}{rl}
 |\left\langle\psi_{m,l,p},\psi_{m,l',p'} \right\rangle|\lesssim_{\g}
\:\:\:\:\:\:\:\:\:\:\:\:\:\: \:\:\:\:\:\:\:\:\:\:\:\:\:\:\:\:\:\:\:\:\:\:\:\:\:\:\:\:\:\:\:\:\:\bigskip\\
|\Phi_{p,p'}''(\Psi_{p,p'}'(\g'(2^{-j})(l-l')))|^{-1/2}\,2^{-m}\,|\phi^{*}(2^{-m}\,\Psi_{p,p'}'(\g'(2^{-j})(l-l')))|\,+\, \textrm{Error}\:,
\end{array}
\eeq
where here $\phi^{*}\in C^N(R)$ with $\textrm{supp}\:\phi^{*}\subseteq\textrm{supp}\:\phi$ while the error term obeys
\beq\label{err}
\textrm{Error}\lesssim \frac{1}{1+|p-p'|+2^{m}\,\g'(2^{-j})\,|l-l'|}\:.
\eeq

\noindent Remark. We used that $\Phi(t)$ defined as above is invertible based on \eqref{fstterm0} and hence we can put ourselves in the setting offered by \eqref{pars}.

We will now analyze the components of the main term in \eqref{interactestim}:

First we want to see under what conditions the main term is nonzero, \textit{i.e.}
\beq\label{difzero}
\phi^{*}(2^{-m}\,\Psi_{p,p'}'(\g'(2^{-j})(l-l')))\not=0\,.
\eeq
Here will be one of the key points where we make use of assumption \eqref{fstterm0}.

Indeed, setting $2^{-m}\,\Psi_{p,p'}'(\g'(2^{-j})(l-l')):=t_0$, we first notice that $t_0\approx 1$.
Next, based on the phase duality, we must have
$\g'(2^{-j})(l-l')=\Phi_{p,p'}'(2^m\,t_0)$. On the other hand using the structure of $\Phi_{p,p'}$ and the mean value theorem we have
$$\Phi_{p,p'}'(t)=R'(\frac{t}{p})-R'(\frac{t}{p'})=t(\frac{1}{p}-\frac{1}{p'})\,R''(c_{pp'})\,$$
with $|t|\approx 2^m$ and $|c_{pp'}|\approx 1$.

 Using now \eqref{fstterm0} we deduce that \eqref{difzero} is equivalent to
 \beq\label{secondmomemtcond}
 2^m\,\g'(2^{-j})\,|l-l'|\approx |p-p'|\:.
 \eeq
 Next we claim that
 \beq\label{below}
 |\Phi_{p,p'}''(\Psi_{p,p'}'(\g'(2^{-j})(l-l')))|\gtrsim_{\g} 2^{-m}\,\g'(2^{-j})\,|l-l'|\;.
 \eeq
 Using that $\Phi_{p,p'},\:\Psi_{p,p'}$ are dual phases, \eqref{below} is a direct consequence of the fact that for $|x|\approx 2^{m}$ we have
 \beq\label{growth}
 |\Phi_{p,p'}''(x)|\gtrsim_{\g} |x|^{-1}\,|\Phi_{p,p'}'(x)|\,
 \eeq
 which further is implied by the hypothesis  \eqref{fstterm0} and \eqref{convdualphase0}.

Thus, taking now $|l-l'|\approx 2^{r}\,|\g'(2^{-j})|^{-1}\,2^{-m}$ (with $r\leq m$) we conclude that in the critical regime (when we have a stationary point in the phase \textit{i.e.} the contribution of the interaction is large) we have
\beq\label{critreg}
|\left\langle\psi_{m,l,p},\psi_{m,l',p'} \right\rangle|\lesssim_{\g} 2^{-r/2}\:.
\eeq

From \eqref{secondmomemtcond} and \eqref{critreg}, we deduce that the main term contribution can be bounded as follows:

\noindent Remark. We will refer
to $D^{Main}_{j,m}$ and $D^{Err}_{j,m}$ as the main and respectively the error term defined in the obvious way by the direct correspondence with \eqref{interactestim}.

$$\|D^{Main}_{j,m}(g,h)\|_2^2\lesssim_{\g} 2^{-m} \,\g'(2^{-j})\sum_{r=0}^{m}\,\sum_{{l,l'=0}\atop{|l-l'|\approx 2^r 2^{-m}}[\g'(2^{-j})]^{-1}}^{[\g'(2^{-j})]^{-1}}\,\sum_{{p,p'=0}\atop{|p-p'|\approx 2^r}}^{2^m} \frac{1}{2^{r/2}}$$
$$| \left\langle g, \varphi_{m,l,p} \right\rangle\ | | \left\langle g, \varphi_{m,l',p'} \right\rangle\ |
| \left\langle h, \varphi_{m,l,p} \right\rangle\ | | \left\langle h, \varphi_{m,l',p'} \right\rangle\ |\lesssim $$
$$2^{-m}\,\g'(2^{-j}) \sum_{r=0}^{m}\,\sum_{{l,l'=0}\atop{|l-l'|\approx 2^r 2^{-m}[\g'(2^{-j})]^{-1}}}^{[\g'(2^{-j})]^{-1}}\,\sum_{{p,p'=0}\atop{|p-p'|\approx 2^r}}^{2^m} \frac{1}{2^{r/2}}
| \left\langle g, \varphi_{m,l,p} \right\rangle\ |^2 | \left\langle h, \varphi_{m,l',p'} \right\rangle\ |^2$$
$$\lesssim 2^{-m}\,\g'(2^{-j})\sum_{r=0}^{m}\,\sum_{l=0}^{[\g'(2^{-j})]^{-1}}\,\sum_{p=0}^{2^m} \frac{1}{2^{r/2}}| \left\langle g, \varphi_{m,l,p} \right\rangle\ |^2\, 2^{r}\,[\g'(2^{-j})]^{-1}\,\|h\|_{\infty}^2$$
$$\lesssim 2^{-m/2}\,\|g\|_{2}^2\,\|h\|_{\infty}^2\:.$$

For the error term we no longer have a condition of the type \eqref{secondmomemtcond} but instead we have the extra decay offered by
\eqref{err}. Thus, proceeding as before we easily deduce that

$$\|D^{Err}_{j,m}(g,h)\|_2^2\lesssim_{\g} 2^{-m} \,\g'(2^{-j})\sum_{r=0}^{m}\,\sum_{{l,l'=0}\atop{|l-l'|\approx 2^r 2^{-m}}[\g'(2^{-j})]^{-1}}^{[\g'(2^{-j})]^{-1}}\,\sum_{{p,p'=0}}^{2^m} \frac{1}{2^{r}}$$
$$| \left\langle g, \varphi_{m,l,p} \right\rangle\ | | \left\langle g, \varphi_{m,l',p'} \right\rangle\ |
| \left\langle h, \varphi_{m,l,p} \right\rangle\ | | \left\langle h, \varphi_{m,l',p'} \right\rangle\ |\lesssim $$
$$\lesssim 2^{-m}\,\g'(2^{-j})\sum_{r=0}^{m}\,\sum_{l=0}^{[\g'(2^{-j})]^{-1}}\,\sum_{p=0}^{2^m} \frac{1}{2^{r}}| \left\langle g, \varphi_{m,l,p} \right\rangle\ |^2\, 2^{r}\,[\g'(2^{-j})]^{-1}\,\|h\|_{\infty}^2$$
$$\lesssim m\,2^{-m}\,\|g\|_{2}^2\,\|h\|_{\infty}^2\:.$$

$\newline$
\textbf{Proof of Proposition 1.}
$\newline$

Our intention here is to show that for $\g'(2^{-j})>2^{-m}$ we have that

$$|\l_{j,m}(f,g,h)|\lesssim_{\g} 2^{-\ep m}\|f\|_2\|g\|_2\|h\|_{\infty}\:$$
$\newline$

Our proof will be split in two cases in accordance with the relative size of  $\g'(2^{-j})$ and $2^{-m}$. (See the Appendix for a a more elaborate discussion of why we have two different cases in this proof.)

$\newline$
\textbf{First regime: $2^{m}\,\g'(2^{-j})$ ``small".}
$\newline$

For this scenario, as described in Section 4, we will use the following discretization  of our operator:

$$B_{j,m}(f,g)(x)=2^{-m}\sum_{l=0}^{[\g'(2^{-j})]^{-1}}\sum_{p=2^m}^{2^{m+1}}\sum_{k=0}^{2^{m}\g'(2^{-j})}
 \left\langle g, \varphi_{m,l,p} \right\rangle\,$$ $$\Q_{m,p}f(l\,2^{m}\g'(2^{-j})+k)\,
 \varphi_{-\log_2 (\g'(2^{-j})),l2^{m}\g'(2^{-j})+k,p2^{-m}[\g'(2^{-j})]^{-1}}(x)\:.$$

Using the dual identity
 $$\l_{j,m}(f,g,h):=\int_{\R}B_{j,m}(f,g)(x)\,\overline{h(x)}\, dx=\int_{\R}f(x)\,\overline{D_{j,m}(g,h)(x)}\, dx\:,$$

the above discretization can be rephrased as follows

$$D_{j,m}(g,h)=2^{-m}\sum_{l=0}^{[\g'(2^{-j})]^{-1}}\sum_{p=2^m}^{2^{m+1}}\sum_{k=0}^{2^{m}\g'(2^{-j})}
\overline{\left\langle g, \varphi_{m,l,p} \right\rangle}\,\left\langle h, \tilde{\varphi}_{m,l,k,p} \right\rangle\,
\psi_{m,l,k,p}\:,$$

where we have set $\tilde{\varphi}_{m,l,k,p}:=\varphi_{-\log_2 (\g'(2^{-j})),l2^{m}\g'(2^{-j})+k,p2^{-m}[\g'(2^{-j})]^{-1}}$ and  $\overline{\hat{\psi}_{m,l,k,p}}(\xi):=
\,e^{-i\,p\,R(\frac{\xi}{\frac{p}{2^m}})}\,e^{i\,(l\,2^{m}\,\g'(2^{-j})+k)\,\xi}\,\phi(\xi)\:.$

Next, as expected, we isolate the interaction term
\beq\label{interact'}
\left\langle\psi_{m,l,k,p},\psi_{m,l',k',p'} \right\rangle=
\int_{\R} e^{i\,[p\,R(\frac{\xi}{\frac{p}{2^m}})-p'\,R(\frac{\xi}{\frac{p'}{2^m}})]}\,
e^{-i\,(\g'(2^{-j})(l-l')\,2^{m}+k-k')\xi}\,\phi(\xi)\,\overline{\phi(\xi)}\,d\xi\:.
\eeq
Further, applying the stationary phase principle (and keeping the same notation as those addressing \eqref{interactestim}) we deduce
\beq\label{interactestim'}
\begin{array}{rl}
 |\left\langle\psi_{m,l,k,p},\psi_{m,l',k',p'} \right\rangle|\lesssim_{\g}\:\:\:\:\:\:\:\:\:\:\:\:\:\:\:\:\:\:\:\:\:\:\:\:\:\:\:\:\:\:\:\:\bigskip\\
|\Psi''(\g'(2^{-j})(l-l'))|^{1/2}\,2^{-m}\,|\phi^{*}(2^{-m}\,\Psi'(\g'(2^{-j})(l-l')))|\,+\, \textrm{Error term}\:.
\end{array}
\eeq
where as before $\phi^{*}\in C^N(R)$ with $\textrm{supp}\:\phi^{*}\subseteq\textrm{supp}\:\phi$.

\noindent Remark. The size of the difference $k-k'$ appearing in the phase of the LHS of
\eqref{interact'} is negligible when applying the stationary phase.

Now, following the same arguments as in the proof of Proposition 2, we have that in the critical regime (\textit{i.e.} main term nonzero)
the condition $|l-l'|\approx 2^{r}\,[\g'(2^{-j})]^{-1}\,2^{-m}$ (with $r\leq m$) implies $|p-p'|\approx 2^{r}$ case in which we further have
\beq\label{critreg'}
|\left\langle\psi_{m,l,k,p},\psi_{m,l',k',p'} \right\rangle|\lesssim_{\g} 2^{-r/2}\:.
\eeq

Using now \eqref{critreg'} we deduce (here we ignore the error term)

$$\|D_{j,m}(g,h)\|_2^2\lesssim_{\g} 2^{-2m} \,\sum_{r=0}^{m}\,\sum_{{l,l'=0}\atop{|l-l'|\approx 2^r 2^{-m}}[\g'(2^{-j})]^{-1}}^{[\g'(2^{-j})]^{-1}}\,\sum_{{p,p'=0}\atop{|p-p'|\approx 2^r}}^{2^m} \,\sum_{k,k'=0}^{2^{m}\g'(2^{-j})}$$
$$\frac{1}{2^{r/2}}\,| \left\langle g, \varphi_{m,l,p} \right\rangle\ | | \left\langle g, \varphi_{m,l',p'} \right\rangle\ || \left\langle h, \tilde{\varphi}_{m,l,k,p} \right\rangle\ | | \left\langle h, \tilde{\varphi}_{m,l',k',p'}\right\rangle\ |$$
$$\lesssim 2^{-2m}\,\sum_{r=0}^{m}\,\sum_{{l,l'=0}\atop{|l-l'|\approx 2^r 2^{-m}[\g'(2^{-j})]^{-1}}}^{[\g'(2^{-j})]^{-1}}\,\sum_{{p,p'=0}\atop{|p-p'|\approx 2^r}}^{2^m} \,\sum_{k,k'=0}^{2^{m}\g'(2^{-j})}$$
$$\frac{1}{2^{r/2}}\,| \left\langle g, \varphi_{m,l,p} \right\rangle\ |^2 | \left\langle h, \tilde{\varphi}_{m,l',k',p'} \right\rangle\ |^2$$
$$\lesssim 2^{-2m}\,\sum_{r=0}^{m}\,\sum_{l=0}^{[\g'(2^{-j})]^{-1}}\,\sum_{p=0}^{2^m} \frac{1}{2^{r/2}}| \left\langle g, \varphi_{m,l,p} \right\rangle\ |^2\, (2^{m}\g'(2^{-j}))^2\,2^{r}\,[\g'(2^{-j})]^{-1}\,\|h\|_{\infty}^2$$
$$\lesssim 2^{-m/2}\,(2^{m}\,\g'(2^{-j}))\,\|g\|_{2}^2\,\|h\|_{\infty}^2\:.$$

Thus we conclude that for $\g'(2^{-j})>2^{-m}$
\beq\label{prop1dot1}
\|D_{j,m}(g,h)\|_2^2\lesssim_{\g} 2^{m/2}\,\g'(2^{-j})\,\|g\|_{2}^2\,\|h\|_{\infty}^2\:.
\eeq
Observe that \eqref{prop1dot1} is very efficient when $\g'(2^{-j})$ close to $2^{-m}$ but very inefficient when $\g'(2^{-j})$ is large (\textit{i.e.} close to $1$). For this second situation we need to choose  a different approach.

$\newline$
\textbf{Second regime: $2^{m}\,\g'(2^{-j})$ ``large".}
$\newline$

 In this case, instead of working with the discretized form of our operator $B_{j,m}$, we follow the argument in \cite{Li}, and choose the compact formulation
\beq\label{compform}
\begin{array}{rl}
&B_{j,m}(f,g)(x)=\\
&2^{-\frac{m}{2}}[\g'(2^{-j})]^{\frac{1}{2}}
 \int_{\R}\int_{\R}\hat{f}(\xi)\,\hat{g}(\eta)\,e^{i (\g'(2^{-j})\xi+\eta)\,x}\,
e^{-i\,2^m\,\eta\,R(\frac{\xi}{\eta})}\,\phi(\xi)\,\phi(\eta)\,d\xi\,d\eta \:.
\end{array}
\eeq

(For a different approach, in which we preserve the discretized model of  $B_{j,m}$, see Section 8.2. in the Appendix. This second approach though, is less general since we need to have good properties of the Weyl sums associated with $\g$ and thus we make the assumption that $\g$ is of polynomial type.)

Returning to \eqref{compform}, our intention is to estimate first the $L^2-$norm of $B_{j,m}$ and hope that the highly oscillatory behavior  of the multiplier will produce some cancelation in the kernel of the operator (this is nothing else then rephrasing the mechanism behind the $TT^{*}$ method). As a consequence we have:
$$\|B_{j,m}(f,g)\|_2^2=2^{-m}\g'(2^{-j})\,\int\int\int\int_{{\R^4}\atop{\g'(2^{-j})\xi+\eta=\g'(2^{-j})\xi_1+\eta_1}}$$
$$\hat{f}(\xi)\,\hat{g}(\eta)\,\overline{\hat{f}(\xi_1)}\,\overline{\hat{g}(\eta_1)}\,
e^{-i\,2^m\,\eta\,R(\frac{\xi}{\eta})}\,\phi(\xi)\,\phi(\eta)\,e^{i\,2^m\,\eta_1\,R(\frac{\xi_1}{\eta_1})}\,
\overline{\phi(\xi_1)}\,\overline{\phi(\eta_1)}\,d\xi\,d\eta\,d\xi_1\,d\eta_1$$
$$=2^{-m}\g'(2^{-j})\,\int\int\int F_{\tau}(\xi)\,G_{\tau}(\eta)\,e^{-i\,2^m\,\varphi_{\tau}(\xi,\eta)}\,d\xi\,d\eta\,d\tau\:,$$

where here we have set $$F_{\tau}(\xi):=\hat{f}\phi(\xi)\,\overline{\hat{f}\phi}(\xi-\tau)\,,\:\:\:\:
G_{\tau}(\eta):=\hat{g}\phi(\eta)\,\overline{\hat{g}\phi}(\eta+\g'(2^{-j})\tau)$$
$$\:\textrm{and}\:\:\:\:\:\varphi_{\tau}(\xi,\eta):=\eta\,R(\frac{\xi}{\eta})-(\eta+\g'(2^{-j})\tau)R(\frac{\xi-\tau}{\eta+\g'(2^{-j})\tau})\:.$$

From the conditions imposed on the curve $\g$ (here is the key point where one makes use of \eqref{convdualphase0})
we deduce that, for large $j$ ($j\,\rightarrow\,\infty$)
\beq\label{hormandphase}
|\partial_{\xi}\partial_{\eta}\varphi_{\tau}(\xi,\eta)|\gtrsim|\tau|\,,
\eeq
relation which implies high oscillatory behavior of the phase and thus enough cancelation for obtaining some decay.
These observation can be captured in the use of a variant of (non) stationary phase principle in two dimensions found in literature under the name of H\"ormander principle (\cite{H}, \cite{S}, Ch. IX); more exactly, we have:
\beq\label{hormand}
|\int\int F_{\tau}(\xi)\,G_{\tau}(\eta)\,e^{-i\,2^m\,\varphi_{\tau}(\xi,\eta)}\,d\xi\,d\eta|\lesssim
\min \{1,|2^m\tau|^{-1/2}\}\,\|F_{\tau}\|_2\,\|G_{\tau}\|_2\:.
\eeq
 Applying now \eqref{hormand} we have
 $$\|B_{j,m}(f,g)\|_2^2$$
 $$\lesssim 2^{-m}\,\g'(2^{-j})\,\sup_{u\in[0,1]}\,
 \{u\,\|f\|_2^2\,\|g\|_2^2 \,+\,(2^m u)^{-1/2}\int_{|\tau|\geq u}\|F_{\tau}\|_2\,\|G_{\tau}\|_2\,d\tau\}$$
 $$\lesssim 2^{-m}\g'(2^{-j})\,\|f\|_2^2\,\|g\|_2^2 \,\sup_{u\in[0,1]}\{u\,+\,u^{-\frac{1}{2}}\,(2^m\,\g'(2^{-j}))^{-\frac{1}{2}} \}\:,$$
 which implies
 \beq\label{l2norm}
 \|B_{j,m}(f,g)\|_2\lesssim (2^{-m}\g'(2^{-j}))^{\frac{1}{2}}\,(2^m\,\g'(2^{-j}))^{-\frac{1}{6}}\,\|f\|_2\,\|g\|_2\:.
 \eeq
Finally, using Cauchy-Schwarz inequality (recall that based on \eqref{localization}, we may reduce the study of $B_{j,m}$ to the case $B_{j,m}=B_{j,m}\,\chi_{[0,\frac{2^{m}}{\g'(2^{-j})}]}$)
\beq\label{csch}
\|B_{j,m}(f,g)\|_1\lesssim(2^m\,\g'(2^{-j}))^{-\frac{1}{6}}\,\|f\|_2\,\|g\|_2\:.
\eeq

Combining now the conclusions of the previous approaches, more exactly, based on \eqref{prop1dot1}
and \eqref{csch} we conclude

$$|\l_{j,m}(f,g,h)|\lesssim_{\g} \min\{ (2^{m/2}\,\g'(2^{-j}))^{\frac{1}{2}},\,(2^m\,\g'(2^{-j}))^{-\frac{1}{6}}\}\,\|f\|_2\|g\|_2\|h\|_{\infty}\:,$$
which implies
$$|\l_{j,m}(f,g,h)|\lesssim_{\g} 2^{-\frac{m}{16}}\,\|f\|_2\|g\|_2\|h\|_{\infty}\:,$$
thus ending the proof of Proposition 1.

\section{Remarks}
$\newline$

In this section we will discuss some possible ways of extending our result as well as some other interesting open problems that naturally arise
from the present paper.

\noindent \textbf{1) The boundedness range of $H_{\Gamma}$.} Our Main Theorem is likely to be phrased for a larger range of indices; thus we do expect a bound of the type $L^p(\R)\times L^q(\R)\:\rightarrow\: L^r(\R)$ with the obvious homogeneity condition $\frac{1}{p}+\frac{1}{q}=\frac{1}{r}$ and for some particular range for $p,q>1$ and $r\geq 1$. Such a bound should be obtained for the local $L^2$ case. It is o interest to investigate the maximal possible range for these indices.

\noindent \textbf{2) Necessary and sufficient conditions for the curve(s) $\g$.} We do not know if it is possible for one to remove conditions \eqref{convdualphase0} and \eqref{convdualphaseinfty} from the definitions of the ``non-flat" curves near $0$ and $\infty$ respectively. Also would be nice if one could replace the asymptotic behavior condition \eqref{asymptotic0},
at which we add \eqref{fstterm0} with some weaker variant resembling \eqref{convtr} (and similarly for the behavior near infinity). For example, a curve of the form $\g(t)= e^{t}$ near infinity is obviously not slowly varying and does not obey \eqref{asymptoticinfty} hence its not in the class $\n\f_{\infty}$.  On the other hand we see no obvious obstruction for why a result as in the Main Theorem shouldn't be transferable to this particular choice of $\g$.

\noindent \textbf{3) More general classes of curves - including the ``flat" case.} Another direction of investigation should concern the limiting situation when our curve becomes very close to being a "line", or with other words infinitely flat. One such example can be provided by $\g(t)=|t|\,\log |t|$. In a situation like this, we do not expect a scale type decay as in the proof of our theorem but rather something closer to the approach of the classical bilinear Hilbert transform where one needs to make use of the concept of tree and organize families of such trees depending on their ``size". A further step is made if one allows $\g$ to also oscillate, like for example
$\g(t)= t\,\sin t$, $\g(t)= t\,e^{\a\,\frac{1}{t}}\,\sin \frac{1}{t}\,,\:\a\in\R$ etc.

Most of these questions can be visualized as part of the program addressing the following
$\newline$

\noindent\textbf{Open problem.} \textit{As before, for $\G=(t,-\gamma(t))$, define the bilinear Hilbert transform $H_{\G}$ along the curve $\G$ as $$H_{\G}: S(\R) \times S(\R)\longmapsto S'(\R)$$
$$H_{\G}(f,g)(x):= \textrm{p.v.}\int_{\R}f(x-t)g(x+\g(t))\frac{dt}{t}\:.$$
What is the largest function space
$$Y\subset C_{H}:=\{\g\in\C(\R\setminus\{0\})\,|\,m(\xi):=p.v.\int_{\R} e^{i\,(t+\g(t))\,\xi}\,\frac{dt}{t}\in L^{\infty}(\R)\}$$
for which the condition $\g\in Y$ assures that $H_{\G}$ extends boundedly from $L^2(\R)\times L^2(\R)$ to $L^1(\R)$?}
$\newline$

\noindent \textbf{4) A curved model for the trilinear Hilbert transform.} Philosophically, we do expect that any general treatment of the above question should distinguish between the flat and non-flat cases and treat them accordingly, combining the method presented here with the deeper and more complicated approach in \cite{LT1}, \cite{LT2}. As a principle, we do notice that the presence of the curvature makes the treatment of this problem easier, as we do obtain a scale-type decay. Accounting for this principle we can introduce the following parallel: the difficult open question regarding the boundedness of the trilinear Hilbert transform
\beq\label{tri}
T(f,g,h)(x):=p.v.\int_{\R}\,f(x+t)\,g(x+2t)\,h(x+3t)\,\frac{dt}{t}
\eeq
say from $L^4(\R)\times L^4(\R)\times L^4(\R)$ to $L^{\frac{4}{3}}(\R)$ should have a significantly easier correspondent "curved" type model as given for example by the question of providing similar bounds for the operator
\beq\label{curvtri}
T_C(f,g,h)(x):=p.v.\int_{\R}\,f(x+t)\,g(x+t^2)\,h(x+t^3)\,\frac{dt}{t}\:.
\eeq
One can further generalize \eqref{curvtri} to trilinear operators integrated over more general families of curves (see also \cite{Li}) or, in another
direction extend these questions to the $n-$linear setting.

\noindent \textbf{5) Another trilinear model problem.} Finally, as another toy model in trying to better understand \eqref{tri} one can ask the problem of providing bounds of the form
$L^p(\R)\times L^q(\R)\times L^r(\R)$ to $L^{s}(\R)$ with $\frac{1}{p}+\frac{1}{q}+\frac{1}{r}=\frac{1}{s}$, $p,q>1$ and $r=\infty$ for the operator given by
\beq\label{bitri}
T(f,g,h)(x):=p.v.\int_{\R}\,f(x+t)\,g(x+2t)\,h(x+t^2)\,\frac{dt}{t}\;.
\eeq
An answer that involves a nontrivial range to this problem would contain some range of exponents from the flat and non-flat bilinear Hilbert transform and hence will have a level of difficulty situated between that of the bilinear Hilbert transform and that of the trilinear one.

$\newline$
\section{Appendix}
$\newline$

\subsection{A comparative discussion}
$\newline$

In this section, for the particular case $\g(t)=t^d$ with $(d\in\N, d\geq 2)$, we briefly make a parallel between
our approach and Li's one.

In \cite{Li} he proves the following
$\newline$

\noindent\textbf{Theorem.} \textit{ If $\g(t)=t^d,\:d\in\N,\:d\geq 2$ then:}

\noindent\textit{For $\g'(2^{-j})>2^{-m}$ we have}
 $$\left\| B_{j,m}(f,g)\right\|_1\lesssim 2^{\frac{(d-1)|j|-m}{8}}\left\|f\right\|_2\left\|g\right\|_2\:.$$
\textit{For $\g'(2^{-j})\leq 2^{-m}$ we have
$$\left\| B_{j,m}(f,g)\right\|_1\lesssim \max\{2^{-\frac{(d-1)|j|-m}{3}},\,2^{-\d\,m}\}\left\|f\right\|_2\left\|g\right\|_2\:,$$
where here $\d>0$ is some absolute constant.}
$\newline$

In Li's approach, the operator $B_{j,m}$ is not discretized but rather embedded in a trilinear expression. The key ingredient used there
is the $\sigma-$uniformity notion introduced in \cite{3}, and inspired by the work of Gowers in \cite{G}.
$\newline$

\noindent \textbf{Observation.} Our intention is to show how the $\sigma-$uniformity concept can be adapted to our discretized setting and naturally bring the same desired conclusion. The entire discussion is striped
as much as possible of technical details exactly for making the analogies and comparisons transparent. However we stress the fact that the ideas presented here can be extended (with some modifications but without deep subtleties) to the general case  of curves $\g$ introduced at the beginning of our paper.
$\newline$

As explained, we will focuss on the boundedness properties of the operator
$$H_{\G}(f,g)(x):= \textrm{p.v.}\int_{\R}f(x-t)g(x+t^d)\frac{dt}{t}\:.$$

With the previous notations, the  $L^2(\R)\times L^2(\R)$ to $L^1(\R)$ boundedness of the
operator $H_{\G}$ follows from our Main Theorem as a consequence of:
$\newline$

\noindent\textbf{Proposition $1'$.}
\textit{For $\g'(2^{-j})>2^{-m}$ we have that}
$$\left\| B_{j,m}(f,g)\right\|_1\lesssim_{\g} 2^{-\frac{m}{16}}\left\|f\right\|_2\left\|g\right\|_2\:.$$

\noindent\textbf{Proposition $2'$.}
\textit{For $\g'(2^{-j})\leq 2^{-m}$ we have that}
$$\|D_{j,m}(g,h)\|_2\lesssim_{\g} 2^{-\frac{m}{4}}\|g\|_2\|h\|_{\infty}\:.$$

\noindent\textbf{Remark}. Notice that our bounds for $\left\| B_{j,m}\right\|_1$ are significantly improved and in particular
are independent of the $j$-parameter. This is a natural consequence of the key ingredient involved in our proof: the careful discretization process that allows the manifestation of the non-vanishing curvature condition \eqref{fstterm0} which further focusses on the second derivative of $\g$ and hence aiming for $2^{-\frac{m}{2}}$-decay. As a consequence
of this process one is able to obtain scale type decay in the $m$ parameter.

In what follows we will only discuss the second proposition.

$\newline$
\subsubsection{\textbf{Proof of Proposition $2'$ (interaction approach).}}
$\newline$

Recall that in the simplified case $\g(t)=t^d,\:d\in\N,\:d\geq 2$ one has

$$D_{j,m}(g,h)=2^{-\frac{m}{2}-\frac{j(d-1)}{2}}\sum_{l=0}^{2^{j(d-1)}}\sum_{p=2^m}^{2^{m+1}}
 \overline{\left\langle g, \varphi_{m,l,p} \right\rangle}\,
 \left\langle h, \varphi_{m,l,p} \right\rangle\,\psi_{m,l,p}(x)\:,$$
 where $\overline{\hat{\psi}_{m,l,p}}(\xi):=\,e^{-i\,\frac{d-1}{d}\, \frac{(2^m \xi)^{\frac{d}{d-1}}}{p^{\frac{1}{d-1}}}}\,e^{i\,d\,l\,2^{m-j(d-1)}\,\xi}\,\phi(\xi)$.

 Further, we notice that
$$\left\langle\psi_{m,l,p},\psi_{m,l',p'} \right\rangle=
\int_{\R} e^{i\,\frac{d-1}{d}\,(\frac{1}{p^{\frac{1}{d-1}}}-\frac{1}{p'^{\frac{1}{d-1}}})\,
(2^m\xi)^{\frac{d}{d-1}}}\,e^{-i\,d\,2^{m-j(d-1)}(l-l')\,\xi}\,\phi(\xi)\,\overline{\phi(\xi)}\,d\xi\:.$$
Consequently, as in \eqref{interactestim}, one has
$$|\left\langle\psi_{m,l,p},\psi_{m,l',p'} \right\rangle|\lesssim_d$$
$$\frac{1}{(1+|p-p'|+2^{m-j(d-1)}|l-l'|)^{\frac{1}{2}}}\:
\phi^{*}\left(c(d)\,\left[\frac{2^{m-j(d-1)}|l-l'|}{|p-p'|}\right]^{d-1}\right)+ \textrm{Error}\:.$$
\noindent Remark. Here $c(d)$ is some real constant depending only on $d$ and $\phi^{*}$ is a smooth,
compactly supported function with $0\notin\textrm{supp}\,\phi^{*}$.

For the main term contribution we have
$$\|D^{Main}_{j,m}(g,h)\|_2^2\lesssim_d 2^{-m-j(d-1)}\sum_{r=0}^{m}\,\sum_{{l,l'=0}\atop{|l-l'|\approx 2^r 2^{(d-1)j-m}}}^{2^{j(d-1)}}\,\sum_{{p,p'=0}\atop{|p-p'|\approx 2^r}}^{2^m} \frac{1}{2^{r/2}}$$
$$| \left\langle g, \varphi_{m,l,p} \right\rangle\ | | \left\langle g, \varphi_{m,l',p'} \right\rangle\ |
| \left\langle h, \varphi_{m,l,p} \right\rangle\ | | \left\langle h, \varphi_{m,l',p'} \right\rangle\ |\lesssim$$
$$\lesssim 2^{-m-j(d-1)} \sum_{r=0}^{m}\,\sum_{{l,l'=0}\atop{|l-l'|\approx 2^r 2^{j(d-1)-m}}}^{2^{j(d-1)}}\,\sum_{{p,p'=0}\atop{|p-p'|\approx 2^r}}^{2^m} \frac{1}{2^{r/2}}
| \left\langle g, \varphi_{m,l,p} \right\rangle\ |^2 \,| \left\langle h, \varphi_{m,l',p'} \right\rangle\ |^2$$
$$ \lesssim2^{-m-j(d-1)}\sum_{r=0}^{m}\,\sum_{l=0}^{2^{j(d-1)}}\,\sum_{p=0}^{2^m} \frac{1}{2^{r/2}}| \left\langle g, \varphi_{m,l,p} \right\rangle\ |^2\, 2^{r+j(d-1)}\,\|h\|_{\infty}^2$$
$$\lesssim 2^{-m/2}\,\|g\|_{2}^2\,\|h\|_{\infty}^2\:.$$
\noindent Remark. The treatment of the error term follows the same lines as in the proof of Proposition 2. We leave these details for the interested reader.

$\newline$

\subsubsection{\textbf{Proof of Proposition $2'$ - translation of Li's approach using $\sigma-$uniformity.}}
$\newline$

Set $\sigma\in (0,1]$ and let $\Q$ be a family of real-valued measurable functions. Also set $I$(=[0,1]) a bounded interval.
$\newline$

\noindent\textbf{Definition.} \textit{A function $f\in L^2(I)$ is $\sigma-$uniform in $\Q$ if}
$$\left|\int_{I} f(\xi)\,e^{-i\,q(\xi)}\,d\xi\right|\leq \sigma \|f\|_{L^2(I)}$$
\textit{for all $q\in\Q$.}
$\newline$

\noindent\textbf{Lemma 0.} \textit{Let $L$ be a bounded sub-linear functional from $L^2(I)$ to $\C$, and $S_{\sigma}$ be the collection of all functions that are $\sigma-$uniform in $\Q$. Set $$U_{\sigma}=\sup_{f\in S_{\sigma}}\frac{|L(f)|}{\|f\|_{L^2(I)}}\:\:\:\&\:\:\:Q=\sup_{q\in\Q}|L(e^{i\,q})|\:.$$
Then for all functions in $L^2(I)$ we have}
$$|L(f)|\leq \max \{U_{\sigma},\,2\sigma^{-1} Q\}\|f\|_{L^2(I)}\:.$$
$\newline$
Now we remind here that for $j(d-1)>m$ we have that
$$ \l_{j,m}(f,g,h):=\int_{\R}B_{j,m}(f,g)(x) \bar{h}(x) dx$$
$$=2^{-\frac{m}{2}-\frac{j(d-1)}{2}}\sum_{l=0}^{2^{j(d-1)}}\sum_{p=2^m}^{2^{m+1}}
 \left\langle g, \varphi_{m,l,p} \right\rangle\,\Q_{m,p}f(l\,2^{m-j(d-1)})\,
 \overline{\left\langle h, \varphi_{m,l,p} \right\rangle}$$
$$=2^{-\frac{m}{2}-\frac{j(d-1)}{2}}\sum_{l=0}^{2^{j(d-1)}}\sum_{p=2^m}^{2^{m+1}}
 \left\langle g, \varphi_{m,l,p} \right\rangle\,
 \left\langle f, \psi_{m,l,p}\right\rangle\ \overline{\left\langle h, \varphi_{m,l,p} \right\rangle}\:,$$
where, as before

 $\overline{\hat{\psi}_{m,l,p}}(\xi):=\,e^{-i\,c_{d}\, \frac{(2^m \xi)^{\frac{d}{d-1}}}{p^{\frac{1}{d-1}}}}\,e^{i\,l\,2^{m-j(d-1)}\,\xi}\,\phi(\xi)\:.$
$\newline$

\begin{p1}
For $j(d-1)\geq m$ we have that
$$| \l_{j,m}(f,g,h)|\lesssim 2^{-\frac{m}{8}}\|f\|_2\|g\|_2\|h\|_{\infty}\:.$$
\end{p1}
Set $$\Q:=\{a \xi^{\frac{d}{d-1}}\,+\,b\xi\,|\,2^{m-100}\leq |a|\leq 2^{m+100}\:\:\&\:\:b\in\R\}\:.$$

\begin{l1}
Let $\hat{f}\phi\in L^2([0,1])$ be $\sigma-$uniform in $\Q$. Then
$$|\l_{j,m}(f,g,h)|\lesssim_d \,\sigma\,\|f\|_2\|g\|_2\|h\|_{\infty}\:.$$
\end{l1}

\begin{l1}
Let $q\in\Q$. Then
$$|\l_{j,m}(\check{(e^{i q}\phi)},g,h)|\lesssim_d 2^{-\frac{m}{4}}\,\|g\|_2\|h\|_{\infty}\:.$$
\end{l1}
Notice that applying now Lemma 0 we have
$$|\l_{j,m}(f,g,h)|\leq_d\max\{\sigma, 2\,\sigma^{-1}2^{-\frac{m}{4}}\}\|f\|_2\|g\|_2\|h\|_{\infty}\:.$$
$$\Rightarrow\:\:\:|\l_{j,m}(f,g,h)|\leq 2^{-\frac{m}{8}}\|f\|_2\|g\|_2\|h\|_{\infty}\:.$$
$\newline$

\noindent\textbf{Proof of Lemma 1.}
$\newline$

We want
$$|\l_{j,m}(f,g,h)|\lesssim \sigma\,\|f\|_2\|g\|_2\|h\|_{\infty}\:,$$
for all $\hat{f}\phi\in L^2([0,1])$ be $\sigma-$uniform in $\Q$
$\Rightarrow\:\:\: | \left\langle f, \psi_{m,l,p} \right\rangle\ |\leq  \sigma\,\|f\|_2\:.$
$$|\l_{j,m}(f,g,h)|\lesssim 2^{-\frac{m}{2}-\frac{j(d-1)}{2}}\sum_{l=0}^{2^{j(d-1)}}\sum_{p=2^m}^{2^{m+1}}
| \left\langle g, \varphi_{m,l,p} \right\rangle\,|
 | \left\langle h, \varphi_{m,l,p} \right\rangle\,| | \left\langle f, \psi_{m,l,p} \right\rangle\ |$$
$$\leq \sigma\,\|f\|_2\, 2^{-\frac{m}{2}-\frac{j(d-1)}{2}} \{\sum_{l=0}^{2^{j(d-1)}}\sum_{p=2^m}^{2^{m+1}}| \left\langle g, \varphi_{m,l,p} \right\rangle\,|^2\}^{\frac{1}{2}} \{\sum_{l=0}^{2^{j(d-1)}}\sum_{p=2^m}^{2^{m+1}}| \left\langle h, \varphi_{m,l,p} \right\rangle\,|^2\}^{\frac{1}{2}}$$
$$\leq \sigma\,\|f\|_2\|g\|_2\|h\|_{\infty}\:.$$
$\newline$

\noindent\textbf{Proof of Lemma 2.}
$\newline$

For $q\in\Q$ we want
$$|\l_{j,m}(\check{(e^{i q}\phi)},g,h)|\lesssim_d 2^{-\frac{m}{4}}\|g\|_2\|h\|_{\infty}\:.$$
Set $q(\xi)=a \xi^{\frac{d}{d-1}}+b\xi$ and let $\hat{f}_{0}=e^{i q}\phi\:.$ Then:
$$\left\langle f_0, \psi_{m,l,p} \right\rangle=\int_{\R}\phi(\xi)\,e^{i\,\left(a-\left(\frac{2^{md}}{p}\right)^{\frac{1}{d-1}}\,c_d\right)
\,\xi^{\frac{d}{d-1}}}\,e^{i\,(b+l 2^{m-j(d-1)})\xi}\,d\xi\:. $$
Applying now stationary phase principle we deduce
\beq\label{statphf0}
|\left\langle f_0, \psi_{m,l,p} \right\rangle|\lesssim_d \frac{1}{(|b+l 2^{m-j(d-1)}|+1)^{\frac{1}{2}}}\:.
\eeq
We may consider $b=0$. Take $N\in\N$ to be later specified.

$\newline$

\noindent\textbf{Case 1}     $\:\:\:l\leq N 2^{j(d-1)-m}\:.$

$$|\l_{j,m}^1(f_0,g,h)|$$
$$:=2^{-\frac{m}{2}-\frac{j(d-1)}{2}}\sum_{l=0}^{N 2^{j(d-1)-m}}\sum_{p=2^m}^{2^{m+1}}
| \left\langle g, \varphi_{m,l,p} \right\rangle\,|
 | \left\langle h, \varphi_{m,l,p} \right\rangle\,| | \left\langle f_0, \psi_{m,l,p} \right\rangle\ |$$
$$\lesssim 2^{-\frac{m}{2}-\frac{j(d-1)}{2}}\|g\|_2\|h\|_{\infty} (N 2^{j(d-1)})^{\frac{1}{2}}= N^{\frac{1}{2}}2^{-\frac{m}{2}}\|g\|_2\|h\|_{\infty}\:.$$
$\newline$

\noindent\textbf{Case 2}     $\:\:\:l\geq N 2^{j(d-1)-m}\:.$

$$|\l_{j,m}^2(f_0,g,h)|$$
$$:=2^{-\frac{m}{2}-\frac{j(d-1)}{2}}\sum_{l=N 2^{j(d-1)-m}}^{2^{j(d-1)}}\sum_{p=2^m}^{2^{m+1}}
| \left\langle g, \varphi_{m,l,p} \right\rangle\,|
| \left\langle h, \varphi_{m,l,p} \right\rangle\,| | \left\langle f_0, \psi_{m,l,p} \right\rangle\ |$$
$$\lesssim 2^{-\frac{m}{2}-\frac{j(d-1)}{2}}\|g\|_2\|h\|_{\infty} 2^{\frac{m+j(d-1)}{2}} N^{-\frac{1}{2}}=N^{-\frac{1}{2}}\|g\|_2\|h\|_{\infty}\:.$$
$\newline$

Combining now cases 1 and 2 we have
$$|\l_{j,m}(f_0,g,h)|\lesssim_d (N^{\frac{1}{2}}2^{-\frac{m}{2}}+\,N^{-\frac{1}{2}})\,\|g\|_2\|h\|_{\infty}\:.$$
Choosing now $N=2^{\frac{m}{2}}$ we obtain that
\beq\label{lcontrol}
|\l_{j,m}(f_0,g,h)|\lesssim_d 2^{-\frac{m}{4}}\,\|g\|_2\|h\|_{\infty}\:.
\eeq

\subsection{Some aspects concerning Proposition 1}
$\newline$

In this section we will discuss on two topics referring to the proof of Proposition 1 (see Section 6):

\begin{itemize}
 \item how can one (at the expense of greater technicalities) still use the discretized version of $B_{j,m}(f,g)$ for the $TT^{*}$ argument in the regime $\g'(2^{-j})>2^{-m}$ with $2^m\,\g'(2^{-j})-$large. This approach though is less general than the one provided in Proposition 1, in the sense that requires some nice number theoretical properties of the discretized curve $\g$, e.g. $\g(t)$ a monomial with $\g(t)=t^d,\,d\in\N,\,d\geq2$;

 \item the necessity of having two different treatments corresponding to the two different regimes in the proof.
\end{itemize}

We will start our discussion addressing \textit{the first topic} in the above list.

 For expository reasons we only address the case $\g(t)=t^2$. Also take for simplicity $j=0$. Then with the discretized version we have that
$$B_{0,m}(f,g)(x)= 2^{-m}\,\sum_{p=2^m}^{2^{m+1}}\sum_{k=0}^{2^m} \left\langle g, \varphi_{m,0,p} \right\rangle
Q_{m,p} f(k)\,\varphi_{0,k,p 2^{-m}}\:.$$
Now when studying the $L^2$-interaction $\left\langle B_{0,m}(f,g), B_{0,m}(f,g) \right\rangle$ first issue is to understand
$$\left\langle\varphi_{0,k,p 2^{-m}}, \varphi_{0,k',p' 2^{-m}} \right\rangle$$
which can be expressed as
$$\int_{\R} e^{-i\,k\,(\xi-2^{-m}p)}\,\hat{\varphi}(\xi-2^{-m}\,p)\,\overline{e^{-i\,k'\,(\xi-2^{-m}p')}\,\hat{\varphi}(\xi-2^{-m}\,p')}\,d\xi$$
$$=e^{-i (k'\,2^{-m}\,p'-k\,2^{-m}\,p)}\sum_{r,r'\in\Z}a_r\,\overline{a_{r'}}\,e^{-i\, r\, 2^{-m}\,p}\,e^{i\, r'\, 2^{-m}\,p'}\,\chi (k-k'-r+r')\:,$$
where here $\chi(\cdot)\in S(\R)$ (Schwartz function) and $\{a_r\}_r\in l^1(\Z)$ are fast decaying complex coefficients.

From last equality we deduce that the main term in the expansion of $\|B_{0,m}(f,g)\|_2^2$ arises when $r=r'=0$ and $k=k'$. Thus we can write

$$\|B_{0,m}(f,g)\|_2^2\approx 2^{-2m}\sum_{p,p'=2^m}^{2^{m+1}}\sum_{k=0}^{2^m} \left\langle g, \varphi_{m,0,p} \right\rangle
\overline{\left\langle g, \varphi_{m,0,p'} \right\rangle}\, e^{i \,k\,2^{-m}\,p}\,e^{-i\,k\,2^{-m}\,p'}$$
$$\int_{\R}\,\int_{\R} \hat{f}(\xi)\,\phi(\xi)\,\overline{\hat{f}(\xi_1)}\,\overline{\phi(\xi_1)}\,e^{-i (\frac{2^{2m}}{2p}\xi^2-\frac{2^{2m}}{2p'}\xi_1^2)}\,e^{i\,k(\xi-\xi_1)}\,d\xi\,d\xi_1$$
Now we notice that the oscillatory term is sensitive to variations $\Delta \xi,\,\Delta \xi_1\approx 2^{-m}$ hence we are invited to express $f\sim \sum_{r=2^m}^{2^{m+1}} \left\langle f, \varphi_{m,0,r} \right\rangle\,\varphi_{m,0,r}$ and thus investigate the interaction
$$\int_{\R}\,\int_{\R} \widehat{\varphi_{m,0,r}}(\xi)\,\phi(\xi)\,\overline{\widehat{\varphi_{m,0,r'}}(\xi_1)}\,\overline{\phi(\xi_1)}\,e^{-i (\frac{2^{2m}}{2p}\xi^2-\frac{2^{2m}}{2p'}\xi_1^2)}\,e^{i\,k(\xi-\xi_1)}\,d\xi\,d\xi_1$$
which after calculations gives
\beq\label{calc}
c_{k,r,p}\,\overline{c_{k,r',p'}}\,2^{-m}\,e^{-i\frac{r^2}{2p}}\,e^{i\frac{kr}{2^m}}\,e^{i\frac{r'^2}{2p'}}\,e^{-i\,\frac{kr'}{2^m}}\:,
\eeq
with
\beq\label{defckrp}
c_{k,r,p}:=\int\hat{\v}(\xi)\,\phi(\frac{\xi+r}{2^m})\,e^{-i\,(\frac{\xi^2}{2p}+\,\xi\,\frac{r}{p})}\,
e^{i\,\frac{k}{2^m}\,\xi}\,d\xi\:.
\eeq
Thus we are led to the following estimate
\beq\label{b0m}
\eeq
$$\|B_{0,m}(f,g)\|_2^2\approx 2^{-3m}\,\sum_{p,p'= 2^m}^{2^{m+1}}\sum_{r,r'=2^{m}}^{2^{m+1}}\,g_p\,\bar{g}_{p'}\,f_{r}\,\bar{f}_{r'}
e^{-i\frac{r^2}{2p}}\,e^{i\frac{{r'}^2}{2p'}}\times$$
$$\sum_{k=0}^{2^m}\,c_{k,r,p}\,\overline{c_{k,r',p'}}\,e^{i\frac{k(r+p)}{2^m}}\,e^{-i\frac{k(r'+p')}{2^m}}\,,$$

where we set $g_p=\left\langle g, \varphi_{m,0,p} \right\rangle$ and $f_r=\left\langle f, \varphi_{m,0,r} \right\rangle$.

Using the shape of the coefficients $c_{k,r,p}$ given in \eqref{defckrp}, one deduces that the main contribution in \eqref{b0m}
is bounded from above by
\beq\label{b0m2}
\eeq
$$\tilde{\M}:=\sup_{{|s|\leq 2^m}\atop{|\xi|,|\xi_1|\lesssim 1}} \left\{ 2^{-2m}\,\sum_{\a,\b=0}^{2^{m-1}}\frac{1}{|\a+\b-s|+1}
\times\right.$$
$$\left.\left|\sum_{r=2^m}^{2^{m+1}}\sum_{p=2^{m}}^{2^{m+1}}
(g_p\,\bar{g}_{p-\a})\,(f_{r}\,\bar{f}_{r-\b})
e^{-i\frac{r^2}{2p}}\,e^{i\frac{(r-\b)^2}{2(p-\a)}}\,e^{-i\,\xi\,\frac{r}{p}}\,e^{i\,\xi_1\,\frac{r-\b}{p-\a}}\right| \right\}\,.$$

Following the reasonings below, one will notice that in order to control $\tilde{\M}$ it is enough to analyze the term
\beq\label{b0m1}
\M:=\sup_{|s|\leq 2^m} 2^{-2m}\,\sum_{\a,\b=0}^{2^{m-1}}\frac{1}{|\a+\b-s|+1}\large|\sum_{r=2^m}^{2^{m+1}}\sum_{p=2^{m}}^{2^{m+1}}
(g_p\,\bar{g}_{p-\a})\,(f_{r}\,\bar{f}_{r-\b})
e^{-i\frac{r^2}{2p}}\,e^{i\frac{(r-\b)^2}{2(p-\a)}}\large| \,.
\eeq
Now using the Van der Corput Lemma (\cite{GK}) we claim that for any $\a,\b \in [1,2^{m-1}]$ we have
\beq\label{VdC}
\frac{1}{2^m}\large|\sum_{p=2^{m}}^{2^{m+1}}\sum_{r=2^m}^{2^{m+1}} a_p\,b_r\,e^{-i\frac{r^2}{2p}}\,e^{i\frac{(r-\b)^2}{2(p-\a)}}\large|
\lesssim\frac{m^{1/4}}{\a^{1/6}}\,\|\{a_p\}\|_2\,\|\{b_r\}_r\|_2\:.
\eeq
Suppose for the moment that this is true; setting then $a_p^{\a}:=g_p\,\bar{g}_{p-\a}$, $b_r^{\b}:=f_{r}\,\bar{f}_{r-\b}$,
$a^{\a}=\{a_p^{\a}\}_{p=2^m}^{2^{m+1}}$ and $b^{\b}=\{b_r^{\b}\}_{r=2^m}^{2^{m+1}}$ we have that
\beq\label{b0m1}
\M\lesssim 2^{-m}\,\sum_{\kappa=0}^{2^{m-1}}\frac{1}{\kappa+1}\sum_{\a=0}^{2^{m-1}}
\min\left\{\frac{1}{2^m}\,\|a^{\a}\|_1\,\|b^{\kappa-\a}\|_1,\:\frac{m^{1/4}}{\a^{1/6}+1}\,\|a^{\a}\|_2\,\|b^{\kappa-\a}\|_2\right\}\:.
\eeq
Using a variational argument we further have
$$\M\lesssim 2^{-m}\,\sum_{\kappa=0}^{2^{m-1}}\frac{1}{\kappa+1}\left\{\frac{1}{2^m}\,
(\sum_{\a=0}^{\a_0}\|a^{\a}\|_1^2)^{\frac{1}{2}}\,(\sum_{\a=0}^{\a_0}\|b^{\kappa-\a}\|_1^2)^{\frac{1}{2}}\right\}$$
$$+\,2^{-m}\,\sum_{\kappa=0}^{2^{m-1}}\frac{1}{\kappa+1}\left\{\frac{m^{1/4}}{{\a_0}^{1/6}+1}
(\sum_{\a=\a_0}^{2^{m-1}}\|a^{\a}\|_2^2)^{\frac{1}{2}}\,(\sum_{\a=\a_0}^{2^{m-1}}\|b^{\kappa-\a}\|_2^2)^{\frac{1}{2}}\right\}$$
and hence
$$\M\lesssim \frac{m}{2^m}\left\{\frac{\a_0}{2^m}+\frac{m^{1/4}}{{\a_0}^{1/6}+1}\right\}\,\|f\|_2^2\,\|g\|_2^2\:.$$
Now properly choosing $\a_0$ (\textit{e.g.} $\a_0=m^{\frac{3}{14}}\,2^{\frac{6m}{7}}$) we deduce that
$$\M\lesssim m^{\frac{17}{14}}\,2^{-m}\,2^{-\frac{m}{7}}\;,$$
and hence we indeed have
$$\left\| B_{0,m}(f,g)\right\|_2\lesssim 2^{-\frac{m}{2}}\, 2^{-\frac{m}{16}}\left\|f\right\|_2\left\|g\right\|_2\:.$$
Let us return now to our claim; notice that \eqref{VdC} is equivalent with
\beq\label{VdcRef}
\frac{1}{2^m}\,\sum_{p,l=2^m}^{2^{m+1}}a_p\,\bar{a}_l \left\{\frac{1}{2^m}\sum_{r=2^m}^{2^{m+1}}
e^{i\,(-\frac{r^2}{2p}+\frac{r^2}{2l}+\frac{(r-\b)^2}{2(p-\a)}-\frac{(r-\b)^2}{2(l-\a)})}\right\}\lesssim \frac{m^{1/2}}{\a^{1/3}}\,\|\{a_p\}\|_2^2\:.
\eeq
Set now
$$ S_{p,l,\a,\b}:=\frac{1}{2^m}\sum_{r=2^m}^{2^{m+1}}
e^{i\,(-\frac{r^2}{2p}+\frac{r^2}{2l}+\frac{(r-\b)^2}{2(p-\a)}-\frac{(r-\b)^2}{2(l-\a)})}=:\frac{1}{2^m}\sum_{r=2^m}^{2^{m+1}} v_r\:.$$
Using Van der Corput lemma, for any $1\leq H\leq 2^m$, we have
$$|S_{p,l,\a,\b}|=|\frac{1}{2^m}\sum_{r=2^m}^{2^{m+1}} v_r|\lesssim \left(\frac{1}{H}\sum_{h=0}^{H-1}
|\frac{1}{2^m}\sum_{r=2^m}^{2^{m+1}} v_{r+h}\,\bar{v}_r|\right)^{\frac{1}{2}}\,+\,\frac{H}{2^m}\;.$$
Let $$V_h:=\frac{1}{2^m}\sum_{r=2^m}^{2^{m+1}} v_{r+h}\,\bar{v}_r\:.$$
Then we have
\beq\label{V}
|V_h|\lesssim \min\left\{\frac{2^{2m}}{h\,\a\,|p-l|},\:1\right\}\:.
\eeq
We have now two cases:

- if $H\leq \frac{2^{2m}}{\a\,|p-l|}$ then $|S_{p,l,\a,\b}|\lesssim \frac{H}{2^m}\leq \frac{2^{m}}{\a\,|p-l|}$\:;

- if $2^{m}\geq H\geq \frac{2^{2m}}{\a\,|p-l|}$ then
$$|S_{p,l,\a,\b}|\lesssim \left(\frac{2^{2m}}{\a |p-l|}\,
\frac{\log H}{H}\right)^{\frac{1}{2}}\,+\,\frac{H}{2^m}\:,$$
and by properly choosing $H$ we have
$$|S_{p,l,\a,\b}|\lesssim m^{\frac{1}{2}}\,\left(\frac{2^{m}}{\a\,|p-l|}\right)^{\frac{1}{3}}\;.$$
Thus we conclude that
\beq\label{S}
|S_{p,l,\a,\b}|\lesssim \min\left\{1,\,m^{\frac{1}{2}}\,\left(\frac{2^{m}}{\a\,|p-l|}\right)^{\frac{1}{3}}\right\}\;,
\eeq
which is enough for proving \eqref{VdcRef}.

For the general case $j<m$ we have that the main term in $\|B_{j,m}\|_2^2$ is given by
$$\bar{\M}:=2^{-3m}\,\sum_{l=0}^{2^{j}}\sum_{k=0}^{2^{m-j}}\sum_{p,p'=2^m}^{2^{m+1}}\sum_{r,r'=2^m}^{2^{m+1}}
\,g_{l,p}\,\bar{g}_{l,p'}\,f_{r}\,\bar{f}_{r'}\,$$
$$\times e^{-i\frac{r^2}{2p}}\,e^{i\frac{r'^2}{2p'}}\,c_{k,r,p,l,j}\,\overline{c_{k,r',p',l,j}}\,e^{i 2^{j-m}\,(p-p')\,(l2^{m-j}+k)}\,e^{i\,(r-r')\,(l\,2^{-j}+k\,2^{-m})}\:,$$
where here
\beq\label{defckrp1}
c_{k,r,p,l,j}:=\int\hat{\v}(\xi)\,\phi(\frac{\xi+r}{2^m})\,e^{-i\,(\frac{\xi^2}{2p}+\,\xi\,\frac{r}{p})}\,
e^{i\,(l\, 2^{-j}+k\,2^{-m})\,\xi}\,d\xi\:.
\eeq
and $g_{l,p}=\left\langle g, \varphi_{m,l,p} \right\rangle$.

Now for $j<\d\,m$ (with $\d\in (0,1)$ small) the main term above can be treated in a similar fashion with the case $j=0$ and thus, making use of the Van der Corput Lemma, one obtains the desired $m-$decay.

We pass now to discussing \textit{the second topic} of our subsection - the natural break of the arguments appearing in the proof of the Proposition 2.

The proof for the regime $2^m\,\g'(2^{-j})$ large is using in the continuous case the H\"ormander principle or in the
discrete setting presented above the Van der Corput lemma; essentially though, in our problem, both are manifestations
of the $TT^{*}$ argument which comes to solve our task of obtaining decay in the $m$ parameter for the expression $\|B_{j,m}\|_1$ by searching for cancelations in the expression $\|B_{j,m}\|_2.$  This method though stops being efficient as $2^m\,\g'(2^{-j})$
approaches $1$, and the explanation is coming below:

Taking for simplicity as before the parabola case (hence $\g'(t)=2t$) we isolate in the main term $\bar{\M}$ above only the summand obtained by taking $l=0$ and $p=p'$. Then this resulting term is in absolute value is at least as large as the absolute value of $$2^{-3m}\,\sum_{p=2^m}^{2^{m+1}}\,\sum_{r,r'=2^m}^{2^{m+1}} |g_p|^2\,f_r\,\bar{f}_{r'}\, e^{i\,\frac{-r^2+r'^2}{2p}}\:.$$
Then as one can easily notice, by taking the supremum restricted to $\|f\|_2\leq 1,\,\|g\|_2\leq 1$ in the above display, the size of the corresponding term is at least $2^{-2m}$.
Thus the best one can reasonable hope is $\|B_{j,m}\|_2\sim 2^{-m}$ and hence because $\textrm{supp}\, B_{j,m}= [0, 2^{m+j}]$ if is to only
use Cauchy-Schwarz we are confined to $\|B_{j,m}\|_1\lesssim 2^{\frac{j-m}{2}}$ which is efficient for $j$ small but bad as $j$ approaches $m$.

On the other hand, if we move to the regime $2^m\,\g'(2^{-j})$ small, using duality, our method was focussed on estimating the size
of $\|D_{j,m} (g,h)\|_2$. There we've made essentially use of the interaction between highly oscillatory terms, but, in the end, by applying \eqref{interactestim'}, we only used the size of this interaction. This is good enough for our aim if we limit ourselves to the above mentioned regime but as soon as we try to use the same argument for the entire case $2^m\,\g'(2^{-j})>1$ this approach fails and we are required to use the full force of \eqref{fouriertr} including the signum and not only on the size of the right hand term in \eqref{fouriertr}.

\end{document}